\documentclass{gtmon_a}
\pdfoutput=1

\proceedingstitle{Proceedings of the School and Conference in Algebraic
Topology (The Vietnam National University, Hanoi, 9-20 August 2004)}
\conferencestart{9 August 2004}
\conferenceend{20 August 2004}
\conferencename{School and Conference in Algebraic Topology}
\conferencelocation{Vietnam National University, Hanoi, Vietnam}

\editor{John Hubbuck}
\givenname{John}
\surname{Hubbuck}

\editor{Nguy\~\ecircumflex{}n H V H\uhorn{}ng}
\givenname{H\uhorn{}ng}
\surname{Nguy\~\ecircumflex{}n}

\editor{Lionel Schwartz}
\givenname{Lionel}
\surname{Schwartz}

\title{M\`ui invariants and Milnor operations}

\author{Masaki Kameko}
\givenname{Masaki}
\surname{Kameko}
\address{Department of Mathematics\\
Faculty of Regional Science\\\newline
Toyama University of International Studies\\
65-1 Higashikuromaki\\
Toyama, 930-1292\\ Japan}
\email{kameko@tuins.ac.jp}
\urladdr{}

\author{Mamoru Mimura}
\givenname{Mamoru}
\surname{Mimura}
\address{Department of Mathematics\\
Faculty of Science\\
Okayama University\\\newline
3-1-1 Tsushima-naka\\
Okayama, 700-8530\\ Japan}
\email{mimura@math.okayama-u.ac.jp}
\urladdr{}

\volumenumber{11}
\issuenumber{}
\publicationyear{2007}
\papernumber{6}
\startpage{107}
\endpage{140}

\doi{}
\MR{}
\Zbl{}

\arxivreference{}



\received{26 Feb 2005}
\revised{12 Sept 2005}
\accepted{14 Sept 2005}
\published{14 November 2007}
\publishedonline{14 November 2007}
\proposed{}
\seconded{}
\corresponding{}
\editor{}
\version{}

\makeatletter
\def\cnewtheorem#1[#2]#3{\newtheorem{#1}{#3}[section]
\expandafter\let\csname c@#1\endcsname\c@thm}

\theoremstyle{plain}
\newtheorem{thm}{Theorem}[section]
\cnewtheorem{prop}[thm]{Proposition}
\cnewtheorem{cor}[thm]{Corollary}
\cnewtheorem{lem}[thm]{Lemma}

\theoremstyle{definition}
\cnewtheorem{defn}[thm]{Definition}
\cnewtheorem{conj}[thm]{Conjecture}
\cnewtheorem{exmp}[thm]{Example}

\theoremstyle{remark}
\cnewtheorem{rem}[thm]{Remark}
\cnewtheorem{rems}[thm]{Remarks}
\cnewtheorem{note}[thm]{Note}
\cnewtheorem{warn}[thm]{Warning}

\makeatother  %  move after \newtheorem block

\begin{document}

\begin{asciiabstract}
We describe Mui invariants in terms of Milnor operations and give a
simple proof for Mui's theorem on rings of invariants of polynomial
tensor exterior algebras with respect to the action of finite general
linear groups. Moreover, we compute some rings of invariants of Weyl
groups of maximal non-toral elementary abelian p-subgroups of exceptional
Lie groups.
\end{asciiabstract}

\begin{htmlabstract}
We describe M&ugrave;i invariants in terms of Milnor operations and
give a simple proof for M&ugrave;i's theorem
on rings of invariants of polynomial tensor exterior algebras
with respect to the action of finite
general linear groups. Moreover, we compute some rings of
invariants of Weyl groups of maximal non-toral elementary
abelian p&ndash;subgroups of exceptional Lie groups.
\end{htmlabstract}

\begin{abstract}
We describe M\`ui invariants in terms of Milnor operations and 
give a simple proof for M\`ui's theorem
on rings of invariants of polynomial tensor exterior algebras 
with respect to the action of finite
general linear groups. Moreover, we compute some rings of 
invariants of Weyl groups of maximal non-toral elementary
abelian $p$--subgroups of exceptional Lie groups. 
\end{abstract}

\maketitle

\section{Introduction} \label{section1}

Let $p$ be a fixed odd prime, $q$ the power of $p$ and 
$\mathbb{F}_q$ the finite field of $q$ elements.
Let 
\[ 
P_n= \mathbb{F}_q[x_{1}, \ldots, x_{n}]
\]
be the polynomial algebra in $n$ variables $x_1, \ldots, x_n$ 
over the finite field $\mathbb{F}_q$.
Let \[
E_{n}^{r}=\Lambda^{r}
(dx_1, \ldots, dx_n)
\]
be the $r^{\mathrm{th}}$ component of the exterior algebra of $dx_1, 
\ldots, dx_n$ over the finite field $\mathbb{F}_q$ and 
let
\[
E_{n}=\bigoplus_{r=0}^{n} E_{n}^{r}
\]
be the exterior algebra of $dx_1, \ldots, dx_n$ over 
$\mathbb{F}_q$.
Let 
\[
P_{n}\otimes E_{n}
\]
 be the polynomial tensor exterior algebra in $n$ variables
$x_{1}, \ldots, x_{n}$ over the finite field $\mathbb{F}_q$.
The general linear group $GL_{n}(\mathbb{F}_q)$ and the special 
linear group $SL_{n}(\mathbb{F}_q)$ act on 
both the polynomial algebra $P_{n}$ and the polynomial tensor 
exterior algebra $P_{n}\otimes E_{n}$.
In \cite{dickson}, the ring of invariants of the polynomial 
algebra is determined by Dickson. 
In \cite{mui}, M\`ui determined the ring of invariants of the 
polynomial tensor exterior algebra 
and described the invariants in terms of determinants.

In the first half of this paper, we describe the invariants in 
terms of Milnor operations and give a simpler proof for M\`ui's 
theorem. 
With the notation in \fullref{section2}, we may state M\`ui's theorems.

\begin{thm}[M\`ui]\label{main1}
The ring of invariants  of the polynomial tensor exterior 
algebra $P_{n}\otimes E_{n}$
with respect to the  action of the special linear group 
$SL_{n}(\mathbb{F}_q)$ is a free 
${P_{n}}^{SL_{n}(\mathbb{F}_q)}$--module with the basis 
$\{1, Q_{I}dx_1\ldots dx_n\}$, where $I$ ranges over $A'_{n}$.
\end{thm}

\begin{thm}[M\`ui] \label{main2}
The ring of invariants of the polynomial tensor exterior algebra 
$P_{n}\otimes E_{n}$
with respect to the  action of the general linear group 
$GL_{n}(\mathbb{F}_q)$ is a free 
${P_{n}}^{GL_{n}(\mathbb{F}_q)}$--module with the basis 
$\{1, e_{n}^{q-2}Q_{I}dx_1\ldots dx_n\}$, where $I$ ranges over 
$A'_{n}$.
\end{thm}

The invariant  $Q_{i_1}\ldots Q_{i_{n-r}}dx_1\ldots dx_n$ in 
\fullref{main1} and \fullref{main2} above is, up to sign,  the 
same as 
the M\`ui invariant $[r\co i_1, \ldots, i_{n-r}]$ in \cite{mui}. 
The first half of this paper has some overlap with M C Crabb's 
work \cite{crabb}. 
However, our point of view on M\`ui invariants seems to be 
different from his.

The second half of this paper is a sequel to the authors' work in 
\cite{kameko} on the invariant theory of Weyl groups 
of maximal non-toral elementary abelian
$p$--subgroup $A$ of simply connected exceptional Lie groups.
For $p$ odd prime, up to conjugation, there are only 6 of them, 
for $p=3$, $A=E_{F_4}^3$, $E_{3E_6}^{4}$, $E_{2E_7}^{4}$, 
$E_{E_8}^{5a}$, $E_{E_8}^{5b}$ and for $p=5$, $A=E_{E_8}^{3}$. 
They and their Weyl groups are described by Andersen et al \cite{andersen}.
We computed the polynomial part of the invariants of Weyl groups 
except for the case $p=3$, $A=E_{E_8}^{5b}$ 
as described by the authors \cite{kameko}. 
In this paper, we compute rings of invariants of polynomial 
tensor exterior algebras 
with respect to the action of Weyl groups except for the case 
$p=3$, $A=E_{E_8}^{5b}$.

In \fullref{section2}, we set up the notation used in the above theorems. 
In \fullref{section3}, we prove \fullref{main1} and \fullref{main2}.
In \fullref{section4}, we state \fullref{main3} and using this theorem, 
we compute rings of invariants of Weyl groups of maximal 
non-toral elementary abelian 
$p$--subgroups of simply connected exceptional Lie groups.
In \fullref{section5}, we prepare for the proof of 
\fullref{main4} and 
\fullref{main3}. 
In \fullref{section6}, we prove \fullref{main4}.
In \fullref{section7}, we prove \fullref{main3}.
In \fullref{sectionapp}, we prove that
the invariant  $Q_{i_1}\ldots Q_{i_{n-r}}dx_1\ldots dx_n$ in 
\fullref{main1} and \fullref{main2} above is, up to sign, equal 
to the M\`ui invariant $[r \co i_1, \ldots, i_{n-r}]$. 

We thank N\,Yagita for informing us that a similar description 
of M\`ui invariants to the above form is also known to him.

\section{Preliminaries} \label{section2}

Let $K_n$ be the field of fractions of $P_n$. 
For a finite set $\{y_1, \ldots, y_r\}$, we denote by 
$\mathbb{F}_q\{y_1, \ldots, y_r\}$ the $\mathbb{F}_q$--vector 
space spanned by
$\{y_1, \ldots, y_r\}$.
Let $GL_{n}(\mathbb{F}_q)$ be the set of $n\times n$ invertible 
matrices with  coefficients in $\mathbb{F}_q$.
We denote by $M_{m,n}(\mathbb{F}_q)$ the set of $m \times n$ 
matrices with coefficients in $\mathbb{F}_q$.
In this paper, we consider the contragredient action of the 
finite general linear group, that is, 
for $g \in GL_n(\mathbb{F}_q)$, we define the action of $g$ on 
$P_n \otimes E_n$ by
\[
g x_i = \sum_{j=1}^{n} a_{i,j}(g^{-1}) x_j, \qua g dx_i = 
\sum_{j=1}^{n} a_{i,j}(g^{-1}) dx_j, 
\]
for $i=1, \ldots, n$ and 
\[
g (x \cdot y) = g(x) \cdot g(y),
\]
for $x$, $y$ in $P_n\otimes E_n$, where $a_{i,j}(g^{-1})$ is the 
entry $(i,j)$ in the matrix $g^{-1}$.
For $x_i$, $dx_i$ in the polynomial tensor exterior algebra $P_n 
\otimes E_n$, 
we define cohomological degrees of $x_i$, $dx_i$ by $\deg 
x_i=2$, $\deg dx_i=1$ for $i=1, \ldots, n$  and we consider 
$P_n\otimes E_n$ as a graded 
$\mathbb{F}_q$--algebra. 

Now, we recall Milnor operations $Q_{j}$ for $j=0, 1, \ldots .$
The exterior algebra 
\[
\Lambda(Q_0, Q_1, Q_2, \ldots)
\]
over $\mathbb{F}_q$, generated by Milnor operations, acts on the 
polynomial tensor exterior algebra 
$P_n\otimes E_n$ as follows;
the Milnor operation $Q_{j}$ is a $P_n$--linear derivation 
\[ 
Q_{j}\co P_{n}\otimes E_{n}^r \to P_{n}\otimes E_{n}^{r-1}
\]
defined by the Cartan formula
\begin{eqnarray*}
Q_{j}(xy)&=&(Q_{j}x) \cdot y +(-1)^{\deg x} x \cdot (Q_{j}y)
\end{eqnarray*}
for $x$, $y$ in $P_n\otimes E_n$ and the unstable conditions
\begin{eqnarray*}
Q_{j}dx_{i}&=&x_{i}^{q^{j}}, \\
Q_{j}x_i&=&0,
\end{eqnarray*}
for $i=1, \ldots, n$, $j\geq 0$.
Thus, the action of Milnor operations $Q_j$ commutes with the 
action of the finite general linear group
$GL_{n}(\mathbb{F}_q)$.

It is also clear that the action of $Q_j$ on $P_n$ trivially 
extends to the quotient field $K_n$ and we may regard $Q_j$ as a 
$K_n$--linear homomorphism \[
Q_j\co K_n\otimes E_n^{r} \to K_n \otimes E_n^{r-1}.
\]

We set up additional notations for handling Milnor operations.

\begin{defn}
For a positive integer $n$, we denote by $S_n$ the set 
\[
\{0, 1, \ldots, n-1\}.
\]
Let $A_{n}$ be the set of subsets of $S_n$.
We denote by $A_{n,r}$ the subset of $A_{n}$ consisting of 
\[
I=\{i_1, \ldots, i_r\}
\]
such that  \[
0\leq i_1<\cdots<i_r< n.
\]
We write $Q_{I}$ for 
\[
Q_{i_1}\ldots Q_{i_r}.
\]
We consider $A_{n,0}$ as the set of empty set $\{\emptyset\}$ 
and define $Q_{\emptyset}$ to be $1$.
It is also convenient for us to define $A'_{n}$ to be the union 
of $A_{n,r}$, where $r$ ranges from $0$ to $n-1$.
\end{defn}

\begin{defn}
Let $I$, $J$ be elements of $A_{n}$.
We define $\mathrm{sign}(I,J)$ as follows.
If $I\cap J\not = \emptyset$, then $\mathrm{sign}(I,J)=0$.
If $I\cap J=\emptyset$ and $I\cup J=\{ k_1, \ldots, k_{r+s}\}$, 
then
$\mathrm{sign}(I,J)$ is the sign of the permutation
\[
\left(\begin{array}{cccccc}
i_1, & \ldots, & i_r, & j_1, & \ldots, & j_s\\
k_1,& \ldots, & k_{r}, & k_{r+1}, & \ldots, & k_{r+s}
\end{array}
\right),
\]
where $I=\{i_1, \ldots, i_r\}$, $J=\{j_1, \ldots, j_s\}$,
$i_1<\cdots<i_r$, $j_1<\cdots<j_s$ and  $k_1<\cdots<k_{r+s}$. 
\end{defn}

The following proposition is immediate from the definition above.

\begin{prop} \label{product}
For $I$, $J$ in $A_{n}$, we have 
\[
Q_{I}Q_{J}=\mathrm{sign}(I,J)Q_{K}=(-1)^{rs}\mathrm{sign}(J,I)Q_K
,
\]
where $K=I\cup J$.
\end{prop}

Finally, using Milnor operations in place of determinants, 
we describe Dickson invariants.
We follow the notation of Wilkerson's paper \cite{wilkerson}.
Let $\Delta_n(X)$, $f_{n}(X)$ be polynomials of $X$ over $P_n$ 
of homogeneous degree $q^{n}$ defined respectively by
\begin{eqnarray*}
\Delta_{n}(X)&=&(-1)^{n}Q_{0}\ldots Q_{n} dx_1\ldots dx_n dX\\
&=&\sum_{i=0}^{n} (-1)^{n-i}(Q_{0}\ldots \widehat{Q}_i \ldots 
Q_{n} dx_1\ldots dx_n) X^{q^{i}},\\
f_n(X)&=&\prod_{x\in \mathbb{F}_q\{x_1, \ldots, x_n\}} (X+x),
\end{eqnarray*}
where the cohomological degrees of $dX$, $X$ are $1$, $2$, and 
$Q_idX=X^{q^i}$, $Q_iX=0$, respectively.

\begin{prop} \label{divisible2}
The polynomial $\Delta_n(X)$ is divisible by the polynomial 
$f_n(X)$ and
\[ e_n(x_1, \ldots, x_n) f_n(X)= \Delta_n(X), \]
where $e_n(x_1, \ldots, x_n)=Q_{0}\ldots Q_{n-1}dx_1\ldots 
dx_n\not =0$.
\end{prop}

\begin{proof}
On the one hand,  both $\Delta_n(X)$ and $f_n(X)$ have all $x 
\in\mathbb{F}_q\{x_1, \ldots, x_n\}$ as roots.
On the other hand, 
the coefficient of $X^{q^n}$ in $\Delta_n(X)$ is 
\[
e_n(x_1, \ldots, x_n)=Q_0\ldots Q_{n-1} dx_1\ldots dx_n
\] and 
$f_n(X)$ is monic.
Since both $\Delta_n(X)$ and $f_n(X)$ have the same homogeneous 
degree $q^n$ as polynomials of $X$, 
we have the required equality.
\end{proof}

Thus, we have the following proposition.

\begin{prop} \label{sum2}
We may express $f_n(X)$ as follows{\rm :}
\[
f_n(X) = \sum_{i=0}^{n} (-1)^{n-i}c_{n,i}(x_1, \ldots, 
x_n)X^{q^{i}},
\]
where
\[
Q_0\ldots \widehat{Q}_i \ldots Q_n dx_1\ldots dx_n = e_n(x_1, 
\ldots, x_n)c_{n,i}(x_1, \ldots, x_n)
\]
and $c_{n,n}(x_1, \ldots, x_n)=1$. 
\end{prop}

The above proposition defines Dickson invariants $c_{n,i}(x_1, 
\ldots, x_n)$ for $i=0, \ldots, n-1$.
When it is convenient and  if there is no risk of confusion, 
we write $e_n$, $c_{n,i}$ for $e_n(x_1, \ldots, x_n)$, 
$c_{n,i}(x_1, \ldots, x_n)$, respectively.

\begin{prop} \label{power2} There holds
\[
c_{n,0}(x_1, \ldots, x_n)=e_{n}(x_1, \ldots, x_n)^{q-1}.
\]
\end{prop}

\begin{proof}
It is clear that 
\[
(Q_{i_1} \ldots Q_{i_n} dx_1\ldots dx_n)^{q}
=Q_{i_1+1}\ldots Q_{i_n+1}dx_1\ldots dx_n\]
and so 
\[
e_{n}(x_1, \ldots, x_n)^{q}=e_{n}(x_1, \ldots, x_n)c_{n,0}(x_1, 
\ldots, x_n). 
\proved
\]
\end{proof}

From the above definitions of $c_{n,i}(x_1, \ldots, x_n)$ and 
$e_n(x_1, \ldots, x_n)$ and from the fact that, 
for $g \in GL_{n}(\mathbb{F}_q)$,
\[
g dx_1\ldots dx_n=\det(g^{-1}) dx_1\ldots dx_n,
\]
it follows that 
\[
g e_{n}(x_1, \ldots, x_n)=\det (g^{-1}) e_n(x_1, \ldots, x_n)
\]
 and that 
\[
g c_{n,i}(x_1, \ldots, x_n)=c_{n,i}(x_1, \ldots, x_n).
\]

Thus, it is clear that ${P_n}^{SL_{n}(\mathbb{F}_q)}$ contains 
$e_{n}$ and
$c_{n,1}, \ldots, c_{n,n-1}$ and that 
${P_n}^{GL_{n}(\mathbb{F}_q)}$ contains 
$c_{n,0}, \ldots, c_{n,n-1}$. Indeed, the following results are 
well-known.
For proofs, we refer the reader to Benson \cite{benson}, Smith
\cite{smith} and Wilkerson \cite{wilkerson}.

\begin{thm}[Dickson] \label{dickson} The ring of invariants 
${P_n}^{SL_{n}(\mathbb{F}_q)}$ is a polynomial algebra generated 
by $c_{n, 1}, \ldots, c_{n,n-1}$ and $e_n$.
\end{thm}

\begin{thm}[Dickson]
The ring of invariants ${P_n}^{GL_{n}(\mathbb{F}_q)}$ is a 
polynomial algebra generated 
by $c_{n,0}, \ldots, c_{n,n-1}$.
\end{thm}

In addition, we need the following proposition.

\begin{prop} \label{projection}
Let 
\[
\pi\co\mathbb{F}_q[x_1, \ldots, x_n]\to \mathbb{F}_q[x_1, 
\ldots, x_{n-1}]
\] 
be the obvious projection.
Then, we have
\[
\pi(e_n(x_1, \ldots, x_n))=0
\]
 and, for $i=1, \ldots, n-1$, 
\[
\pi(c_{n,i}(x_1, \ldots, x_n))=c_{n-1, i-1}(x_1,\ldots, 
x_{n-1})^{q}.
\]
\end{prop}

\begin{proof}
It is clear that $e_n(x_1, \ldots, x_n)$ is divisible by $x_n$, 
so we have
\[
\pi(e_n(x_1, \ldots, x_n))=0.
\]
On the one hand, we have
\begin{eqnarray*}
f_n(X) & = & \prod_{\alpha\in \mathbb{F}_q} \qua \prod_{ x\in 
\mathbb{F}_q\{x_1, \ldots, x_{n-1} \}} (X+\alpha x_n + x)\\
& = & \prod_{\alpha\in \mathbb{F}_q}f_{n-1}(X+\alpha x_n)\\
&=&\prod_{\alpha \in \mathbb{F}_q} (f_{n-1}(X)+\alpha 
f_{n-1}(x_n))\\
&=& f_{n-1}(X)^{q}-f_{n-1}(X)f_{n-1}(x_n)^{q-1}.
\end{eqnarray*}
On the other hand, since $f_{n}(x_n)$ is divisible by $x_{n}$, 
we have
\[
\pi(f_{n}(X))=f_{n-1}(X)^{q}.
\]
Comparing the coefficients of $X^{q^{i}}$, we have, for $i=1, 
\ldots, n-1$,
\[
\pi(c_{n,i}(x_1, \ldots, x_n))=c_{n-1,i-1}(x_1, \ldots, 
x_{n-1})^{q}. 
\proved
\]
\end{proof}

\section
[Proof of Theorems 1.1 and 1.2]
{Proof of \fullref{main1} and \fullref{main2}} \label{section3}

In the case $n=1$, the invariants are obvious. In the case 
$r=0$, the invariants are 
calculated by Dickson. So, throughout the rest of this section, 
we assume $n\geq 2$ and $r>0$.
To prove \fullref{main1} and \fullref{main2}, it suffices to 
prove the following theorems.

\begin{thm}\label{special}
The submodule $(P_{n}\otimes E_n^r)^{SL_{n}(\mathbb{F}_q)}$
is a free ${P_n}^{SL_{n}(\mathbb{F}_q)}$--module with 
the basis $\{ Q_{I}dx_1\ldots dx_n \}$, where $I$ ranges over 
$A_{n,n-r}$.
\end{thm}

\begin{thm}\label{general}
The submodule 
$(P_n\otimes E_n^r)^{GL_{n}(\mathbb{F}_q)}$ is a free 
${P_n}^{GL_n(\mathbb{F}_q)}$--module with the basis 
$\{ e_{n}^{q-2}Q_{I}dx_1\ldots dx_n \}$,
where $I$ ranges over $A_{n,n-r}$.
\end{thm}

To begin with, we prove the following proposition.

\begin{prop}\label{basis}
The elements $Q_{I}dx_1\ldots dx_n$ form a basis for $K_n 
\otimes E_n^{r}$, where $I$ ranges over $A_{n, n-r}$.
\end{prop}

\begin{proof}
Firstly, we show the linear independence of $Q_Idx_1\ldots 
dx_n$.
Suppose that 
\[
a = \sum_{I \in A_{n, n-r}} a_{I} Q_{I} dx_1\ldots dx_n,
\]
where $a_{I}\in K_n$.
For each $I\in A_{n, n-r}$, let $J=S_{n}\backslash I$. 
It is clear that $J \cap I'\not = \emptyset$ if $I'\not =I$ in 
$A_{n,n-r}$. Hence, 
we have
$\mathrm{sign}(J,I)\not =0$ and $\mathrm{sign}(J,I')=0$
for $I'\not =I \in A_{n,n-r}$.
By \fullref{product}, we have
\[ {Q}_{J}a = \mathrm{sign}(J,I)a_{I}Q_0\ldots Q_{n-1}dx_1\ldots 
dx_n
=\mathrm{sign}(J,I)a_{I}e_n.\]
Thus, if $a=0$, then $a_{I}=0$. 
Therefore, the terms $Q_{I}dx_1\ldots dx_n$ are linearly independent in 
$K_n \otimes E_{n}^r$.

On the other hand, since 
\[
\dim_{K_n} K_n\otimes E_n^r=\left( \begin{array}{c} n \\ r 
\end{array}\right)
\]
is equal to the number of elements in $A_{n, n-r}$, we see,
for dimensional reasons, that $Q_{I}dx_1\ldots dx_n $'s form a 
basis for $K_n\otimes E_n^r$.
\end{proof}

\begin{lem}\label{lem}
Let $h_{I}$ be polynomials over $\mathbb{F}_q$ in $(n-1)$ 
variables, where $I$ ranges over $A_{n,n-1}$.
Suppose that
\[
a_{0}=\sum_{I\in A_{n,n-1}} h_{I}(c_{n,n-1}, \ldots, c_{n, 
1}){e_{n}}^{-1}
Q_{I} dx_{1}\ldots dx_{n}
\]
is in $P_n\otimes E_n^1$.
Then $h_{I}=0$ for each $I\in A_{n,n-1}$.
\end{lem}

\begin{proof}[Proof of \fullref{special}]

Suppose that $a$ is an element in $P_n\otimes E_n^r$ and that
$a$ is ${SL_{n}(\mathbb{F}_q)}$--invariant.
By \fullref{basis},  the elements $Q_{I}dx_1\ldots dx_n$ form 
a basis for $K_n\otimes E_n^r$. Hence, there exist $a_I$ in 
$K_n$  such that
\[ a=\sum_{I \in A_{n, n-r}} a_{I}Q_{I}dx_1\ldots dx_n. \]
For  $I\in A_{n,n-r}$, let $J=S_n\backslash I$. 
Then, $Q_{J}a$ is in $P_n$.
As in the proof of \fullref{basis},  we have
\[ Q_{J}a=\mathrm{sign}(J,I)a_{I}e_n. \]
Therefore, there are polynomials $f_{I,k}$ over $\mathbb{F}_q$ 
in $(n-1)$ variables such that
\[
a_{I}=\sum_{k\geq 0} f_{I, k}(c_{n,n-1}, \ldots, c_{n,1}) 
e_{n}^{k-1}.
\]
Thus, we have
\[
a= \sum_{I\in A_{n, n-r}} \sum_{k\geq 0} f_{I, k}(c_{n,n-1}, 
\ldots, c_{n,1}) 
e_{n}^{k-1} Q_{I}dx_1\ldots dx_n.
\]
It remains to show that $f_{I,0}=0$ for each $I \in A_{n,n-r}$.

Let
\[
a_0=a-\sum_{I\in A_{n,n-r}} \sum_{k\geq 1} f_{I,k}(c_{n,n-1}, 
\ldots, c_{n,1})e_n^{k-1}Q_Idx_1\ldots dx_n.
\]
Then, we have that 
\[
a_{0} = \sum_{I\in A_{n,n-r}} f_{I,0}(c_{n,n-1}, \ldots, c_{n, 
1})e_n^{-1}Q_{I}dx_1\ldots dx_n
\]
and that $a_0$ is also in $P_n \otimes E_n^r$.

For $J \in A_{n, r-1}$, the element 
$Q_{J}a_{0}$ is in $P_n \otimes E_n^1$. By \fullref{product}, we 
have
\[
Q_{J}a_{0}=\sum_{K\in A_{n,n-1}, I=K\backslash J, I \in A_{n, 
n-r}}
\mathrm{sign}(J,I) f_{I, 0}(c_{n, n-1}, \ldots, c_{n, 1}) Q_{K} 
dx_1\ldots dx_n.
\]
Hence, by \fullref{lem}, we have $\mathrm{sign}(J,I)f_{I,0}=0$.
For each $I$ in $A_{n,n-r}$, there exists $J\in A_{n,r-1}$ such 
that $\mathrm{sign}(J,I)\not=0$.
Therefore, we have $f_{I,0}=0$ for each $I$. This completes the 
proof.
\end{proof}

\begin{proof}[Proof of \fullref{general}]
As in the proof of \fullref{special}, if 
$a \in P_n\otimes E_n^r$ is ${GL_n(\mathbb{F}_q)}$--invariant, 
the element $a$ can be expressed in the form
\[
a=\sum_{I\in A_{n, n-r}} \sum_{k\geq 1} f_{I,k}(c_{n, n-1}, 
\ldots, c_{n,1})e_{n}^{k-1}Q_{I}dx_1\ldots dx_n.
\]
For $g\in GL_n(\mathbb{F}_q)$, we have
\[
ga = \sum_{I\in A_{n, n-r}} \sum_{k\geq 1} 
\det (g^{-1})^{k} f_{I,k}(c_{n, n-1}, \ldots, 
c_{n,1})e_{n}^{k-1}Q_{I}dx_1\ldots dx_n.
\]
Therefore, $a$ is ${GL_{n}(\mathbb{F}_q)}$--invariant if and only 
if 
$f_{I,k}=0$ for $k\not \equiv 0$ modulo $q-1$. Hence, we have
\[a =\sum_{I \in A_{n, n-r}} \sum_{m\geq 0}  f_{I, 
m(q-1)+(q-1)}(c_{n,n-1}, \ldots, c_{n,1}) e_{n}^{m(q-1)} 
e_{n}^{q-2}Q_{I}dx_1\ldots dx_n.
\]
Since $e_{n}^{q-1}=c_{n,0}$, we may write
\[
a =\sum_{I\in A_{n, n-r}} f'_{I}(c_{n, n-1}, \ldots, 
c_{n,1},c_{n,0})e_{n}^{q-2}Q_{I}dx_1\ldots dx_n,
\]
where \[f'_{I}(c_{n,n-1}, \ldots, c_{n,0})
=\sum_{m\geq 0} f_{I, m(q-1)+(q-1)}(c_{n,n-1}, \ldots, 
c_{n,1}){c_{n,0}}^{m}.\]
This completes the proof.
\end{proof}

\begin{proof}[Proof of \fullref{lem}]
For the sake of notational simplicity, let $I_i=S_n\backslash 
\{i\}$ and we write $h_i$ for $h_{I_i}$.
Since $a$ is in $P_n \otimes E_n^1$, 
there are $\varphi_1, \ldots, \varphi_n$ in $P_n$ such that
\[
a_0 = \varphi_1dx_1+\cdots+\varphi_n dx_n.
\]
The coefficient $\varphi_n$ of $dx_{n}$ is given by
\begin{eqnarray*}
&&\sum_{i=0}^{n-1} 
h_{i}(c_{n, n-1}, \ldots, c_{n,1}) e_{n}^{-1}Q_{I_i}dx_1\ldots 
dx_{n-1}\\
&=& \left\{\sum_{i=0}^{n-1} h_{i}(c_{n, n-1}, \ldots, c_{n,1}) 
c_{n-1, i}\right\}
e_{n}^{-1}e_{n-1}.
\end{eqnarray*}
Hence, we have 
\[
\left\{\sum_{i=0}^{n-1} h_{i}(c_{n, n-1}, \ldots, 
c_{n,1})c_{n-1,i}\right\}
e_{n-1}=e_{n}\varphi_n.
\]
By \fullref{projection}, the obvious projection 
\[\pi\co\mathbb{F}_q[x_1, \ldots, x_n] \to \mathbb{F}_q[x_1, 
\ldots, x_{n-1}]
\] maps
$e_{n}$, $c_{n,i}$ to $0$, $c_{n-1,i-1}^{q}$, respectively. 
So, we have
\[
\sum_{i=0}^{n-1} h_{i}(c_{n-1, n-2}^{q}, \ldots, c_{n-1, 0}^{q}) 
c_{n-1,i}=0.
\]
Since $c_{n-1,i}$ $(i=0, \ldots, n-2)$ are algebraically 
independent in $\mathbb{F}_q[x_1, \ldots, x_{n-1}]$ and 
since $c_{n-1, n-1}=1$, writing $y_i$ for $c_{n-1,i}$, we have 
the following equation: 
\begin{equation}
\label{eqn1}
h_{n-1}(y_{n-2}^{q}, \ldots, y_{0}^{q})+\sum_{i=0}^{n-2} 
h_{i}(y_{n-2}^{q}, \ldots, y_{0}^{q})y_{i}=0.
\end{equation}
Applying partial derivatives $\partial/\partial y_i$,  we have
\begin{equation}
h_{i}(y_{n-2}^{q}, \ldots, y_{0}^{q})=0
\end{equation}
for $i=0, \ldots, n-2$.
Hence, $h_{i}(y_{n-2}^{q}, \ldots, y_{0}^{q})=0$ for $i=0, 
\ldots, n-2$. 
Substituting these to the previous equation \eqref{eqn1}, we also have 
$h_{n-1}(y_{n-2}^{q}, \ldots, y_{0}^{q})=0$.
Since we deal with polynomials over the finite field 
$\mathbb{F}_q$, we have
\[
h_{i}(y_{n-2}^{q}, \ldots, y_{0}^{q})=h_{i}(y_{n-2}, \ldots, 
y_{0})^{q}
\]
for $i=0, \ldots, n-1$.
Therefore, we have
$h_{i}(y_{n-2}, \ldots, y_{0})=0$ for $i=0, \ldots, n-1$.
Since $y_0, \ldots, y_{n-2}$ are algebraically independent, we 
have
\[
h_i=0
\] as polynomials over $\mathbb{F}_q$ in $(n-1)$ variables for 
$i=0, \ldots, n-1$.
\end{proof}

\section{Invariants of some Weyl groups} \label{section4}

In this section, we consider the invariant theory of polynomial 
tensor exterior algebras.
In what follows, we assume that $n\geq 2$.
To state \fullref{main3}, which is our main theorem on the 
invariant theory, we need some notation.
Let 
\[
P_{n-1}=\mathbb{F}_{q}[x_2, \ldots, x_n]
\]
be the subalgebra of $P_n$ generated by $x_2, \ldots, x_n$ and
let
$
E_{n-1}
$
be the subalgebra of $E_n$ generated by $dx_2, \ldots, dx_n$.
Let $G_1$ be a subgroup of $SL_{n-1}(\mathbb{F}_q)$ which acts 
on $P_{n-1}$ and $P_{n-1}\otimes E_{n-1}$ both.
Let $G$ be a subgroup of $SL_{n}(\mathbb{F}_q)$ consisting of 
the following matrices:
\[
\left( \begin{array}{c|c}
1 & m \\ \hline
0 &   g_1
\end{array}\right),
\]
where $g_1\in G_1$ and $m\in M_{1,n-1}(\mathbb{F}_q)$. 
Obviously the group $G$ acts on $P_{n}$ and $P_{n}\otimes E_{n}$.
Finally, let 
\[
\displaystyle \mathcal{O}_{n-1}(x_i)=\prod_{x \in 
\mathbb{F}_{q}\{ x_2, \ldots, x_n\}} (x_i+x).
\]

\begin{thm}\label{main4}
Suppose that the ring of invariants 
$
P_{n-1}^{G_1}
$
is a polynomial algebra generated by homogeneous polynomials 
$f_2, \ldots, f_n$ in $(n-1)$ variables 
$x_2, \ldots, x_n$. Then, the ring of invariants 
$
P_{n}^{G}
$
is also a polynomial algebra generated by 
$$\mathcal{O}_{n-1}(x_1), f_{2}, \ldots, f_n.$$
\end{thm}

This theorem is a particular case of a theorem of Kameko and Mimura 
\cite[Theorem 2.5]{kameko}.
We use this theorem to compute the polynomial part of invariants 
$P_n^{G}$ which appear in our main theorem, \fullref{main3}.
So, \fullref{main3} below works effectively together with 
\fullref{main4}.

\begin{thm}\label{main3}
Suppose that the ring of invariants
$
P_{n-1}^{G_1}
$
 is a polynomial algebra generated by homogeneous polynomials 
$f_2, \ldots, f_{n}$
in $(n-1)$ variables $x_2, \ldots, x_n$ and
suppose that the ring of invariants 
$
(P_{n-1} \otimes E_{n-1})^{G_1}
$
is a free $P_{n-1}^{G_1}$--module with a basis 
$
\{ v_i \},
$
where $i=1, \ldots, 2^{n-1}$.
Then, the ring of invariants 
$
(P_{n}\otimes E_{n})^{G}
$
is a free $P_{n}^{G}$--module with the basis 
$
\{ v_{i}, Q_{I} dx_1\ldots dx_n  \},
$
where $i=1, \ldots, 2^{n-1}$ and  $I$ ranges over $A_{n-1}$.
\end{thm}

We prove \fullref{main4} and \fullref{main3} in \fullref{section5}, 
\fullref{section6} and  \fullref{section7}.

As an application of \fullref{main4} and \fullref{main3}, 
we compute rings of invariants of
the mod--$p$ cohomology of the classifying spaces of maximal 
non-toral elementary abelian $p$--subgroups
of simply connected  exceptional Lie groups with respect to the 
Weyl group action.

It is well-known that for an odd prime $p$, a simply connected  
exceptional Lie group $G$ does not have non-toral elementary 
abelian $p$--subgroups 
except for the cases $p=5$,  $G=E_8$, and $p=3$, $G=F_4$, $E_6$, 
$E_7$, $E_8$
(see \cite{andersen} and \cite{greiss}).
Andersen, Grodal, M{\o}ller and Viruel \cite{andersen} described 
the Weyl groups of maximal 
 non-toral elementary abelian $p$--subgroups
as well as their action on the underlying elementary abelian 
$p$--subgroup explicitly for $p=3$, $G=E_6$, $E_7$, $E_8$.
Up to conjugate, there are only 6 maximal non-toral elementary 
abelian $p$--subgroups of simply connected  exceptional Lie 
groups.
For $p=5$, $G=E_8$ and for $p=3$, $G=F_4$, $E_6$, $E_7$, 
there is one maximal non-toral elementary abelian $p$--subgroup 
for each $G$. We call them
$E_{E_8}^3$, $E_{F_4}^3$, $E_{3E_6}^4$, $E_{2E_7}^4$, following 
the notation in \cite{andersen}.
For $p=3$, $G=E_8$, there are two maximal non-toral elementary 
abelian $p$--subgroups, say $E_{E_8}^{5a}$ and $E_{E_8}^{5b}$,
where the superscript indicates the rank of elementary abelian 
$p$--subgroup.
For a detailed  account on non-toral elementary abelian 
$p$--subgroups, 
we refer the reader to Andersen et al 
\cite[Section 8]{andersen}, and its references.

Let $A$ be an elementary abelian $p$--subgroup of a compact Lie 
group $G$. 
Suppose that $A$ is of rank $n$.
We denote by $W(A)$ the Weyl group 
of $A$. Choosing a basis, say $\{ a_i\}$, for $A$, we consider 
the Weyl group $W(A)$ as 
a subgroup of the finite general linear group 
$GL_n(\mathbb{F}_p)$.
We write $H^{*}BA$ for the mod--$p$ cohomology of the classifying 
space $BA$.
The Hurewicz homomorphism $h\co A=\pi_1(BA) \to 
H_1(BA;\mathbb{F}_p)$ 
is an isomorphism.  We denote by $\{ dt_i\}$ the dual basis of 
$\{ h(a_i) \}$, so that $dt_i$ is the dual of $h(a_i)$ with 
respect to
the basis $\{h(a_i)\}$ of $H_1(BA;\mathbb{F}_p)$ for $i=1, 
\ldots, n$. Let
$\beta\co H^1 BA \to H^2BA$ be the Bockstein homomorphism.
Then, the mod--$p$ cohomology of $BA$ is a polynomial tensor 
exterior algebra
\[
H^{*}BA = \mathbb{F}_p[t_1, \ldots, t_n]  \otimes 
\Lambda(dt_1, \ldots, dt_n),
\]
where $\deg t_i=2$, $\deg dt_i=1$ and $t_i=\beta(dt_i)$ for 
$i=1, \ldots, n$. 
We denote by $\Gamma H^{*}BA$ the polynomial part of
$H^{*}BA$, that is, 
\[
\Gamma H^{*}BA=\mathbb{F}_p[t_1, \ldots, t_n].
\]
The action of the Weyl group $W(A)$ on $A=\pi_1(BA)$, given by
\[
g a_i = \sum_{j} a_{j,i}(g) a_j,
\]
where $\{a_i\}$ is the fixed basis of $A$, induces the action of 
$W(A)$ on $H^{*}BA$, which is 
given by
\[
g t_i= \sum_{j} a_{i,j}(g^{-1}) t_j, \qua g dt_i =\sum_{j} 
a_{i,j}(g^{-1}) dt_j,
\]
for $i=1, \ldots, n$. 

Now, we compute
\[
(H^{*}BA)^{W(A)}
\]
for $A=E_{E_8}^{3}$, $E_{F_4}^{3}$, $E_{3E_6}^4$, $E_{2E_7}^4$, 
$E_{E_8}^{5a}$ using \fullref{main1}, \fullref{main4} 
and \fullref{main3}.  

\begin{prop}\label{ring} 
For the above elementary abelian $p$--subgroup $A$, the ring of 
invariants 
$(H^{*}BA)^{W(A)}$ is given as follows{\rm :}
\begin{itemize}
\item[{\rm (1)}] For $p=5$, $G=E_8$, $\displaystyle 
A=E_{E_8}^3$, 
$
(H^*BA)^{W(A)}$
is given by
\[
\mathbb{F}_5[x_{62}, x_{200}, x_{240}] \otimes \mathbb{F}_5 \{1, 
Q_I u_3\}, 
\]
where $x_{62}=e_{3}(t_1, t_2, t_3)$, $x_{200}=c_{3,2}(t_1, t_2, 
t_3)$, $x_{240}=c_{3,1}(t_1, t_2, t_3)$, 
$u_3=dt_1dt_2dt_3$ and $I$ ranges over $A_3'${\rm ;}
\item[{\rm (2)}] For $p=3$, $G=F_4$, $A=E_{F_4}^3$,
$(H^*BA)^{W(A)}$ is given by
\[
\mathbb{F}_3[x_{26}, x_{36}, x_{48}] \otimes \mathbb{F}_3\{
1, Q_{I}u_3\},
\]
where $x_{26}=e_{3}(t_1, t_2, t_3)$, $x_{36}=c_{3,2}(t_1, t_2, 
t_3)$, $x_{48}=c_{3,1}(t_1, t_2, t_3)$, 
$u_3=dt_1dt_2dt_3$ and $I$ ranges over $A_3'${\rm ;}
\item[{\rm (3)}] For $p=3$, $G=E_6$, $A=E_{3E_6}^{4}$, 
$
(H^*BA)^{W(A)}$
is given by
\[
\mathbb{F}_3[x_{26}, x_{36}, x_{48}, x_{54}]
\otimes \mathbb{F}_3\{1, Q_{I}u_3, Q_{J}u_4\},
\]
where $x_{26}=e_{3}(t_2, t_3, t_4)$, $x_{36}=c_{3,2}(t_2, t_3, 
t_4)$, 
$x_{48}=c_{3,1}(t_2, t_3, t_4)$, 
\[
x_{54}=\prod_{t\in \mathbb{F}_3\{t_2, t_3, t_4\}} (t_1+t), 
\]
$u_3=dt_2dt_3dt_4$, $u_4=dt_1dt_2dt_3dt_4$, $I$ ranges over 
$A_{3}'$ and $J$ ranges over $A_{3}${\rm ;}
\item[{\rm (4)}] For $p=3$, $G=E_7$, $A=E_{2E_7}^{4}$, 
$
(H^*BA)^{W(A)}$
is given by
\[
\mathbb{F}_3[x_{26}, x_{36}, x_{48}, x_{108}]
\otimes \mathbb{F}_3
\{1, Q_I u_3, x_{54} Q_{J} u_4\},
\]
where $x_{26}=e_{3}(t_2, t_3, t_4)$, $x_{36}=c_{3,2}(t_2, t_3, 
t_4)$, 
$x_{48}=c_{3,1}(t_2, t_3, t_4)$, $x_{108}=x_{54}^2$, 
\[
\displaystyle x_{54}=\prod_{t \in \mathbb{F}_3\{ t_2, t_3, 
t_4\}}(t_1+t),
\]
$u_3=dt_2dt_3dt_4$, $u_4=dt_1dt_2dt_3dt_4$, $I$ ranges over 
$A_{3}'$ and $J$ ranges over $A_{3}${\rm ;}
\item[{\rm (5)}] For $p=3$, $G=E_8$, $A=E_{E_8}^{5a}$, 
$
(H^*BA)^{W(A)}$
is given by
\[
\mathbb{F}_3[x_{4}, x_{26}, x_{36}, x_{48}, x_{324}] \otimes
\mathbb{F}_3\{ 1, Q_{I}u_3, x_2u_1, (Q_{I}u_3)x_2u_1, 
x_2x_{162}Q_{J}u_5\},
\]
where $x_4=x_2^2$, $x_{26}=e_{3}(t_2, t_3, t_4)$, 
$x_{36}=c_{3,2}(t_2, t_3, t_4)$, 
$x_{48}=c_{3,1}(t_2, t_3, t_4)$, $x_{324}=x_{162}^2$, $x_2=t_5$, 
\[
x_{162}=\prod_{t\in \mathbb{F}_3\{t_2, t_3, t_4, t_5\}} (t_1+t),
\]
$u_1=dt_5$, $u_3=dt_2dt_3dt_4$, $u_5=dt_1dt_2dt_3dt_4dt_5$, 
$I$ ranges over $A_{3}'$ and  $J$ ranges over $A_{4}$, 
\end{itemize}
where the subscripts of $u$ and $x$  indicate their 
cohomological degrees.
\end{prop}

\begin{proof}
(1), (2)\qua
In the case $p=5$, $G=E_8$, $A=E_{E_8}^3$ and 
in the case $p=3$, $G=F_4$, $A=E_{F_4}^3$, the Weyl group is 
$SL_3(\mathbb{F}_p)$. Therefore, 
it is the case of M\`ui invariants and it is immediate from 
\fullref{main1}. 

(3)\qua In the case $p=3$, $G=E_6$, $A=E_{E_6}^4$, the Weyl group 
$W(A)$ is the subgroup of $SL_{4}(\mathbb{F}_3)$
consisting of the following matrices:
\[
\left( \begin{array}{c|c} 
1 & m \\ \hline
0 & g_1 
\end{array}
\right),
\]
where $m\in M_{1,3}(\mathbb{F}_3)$, $g_1 \in 
SL_{3}(\mathbb{F}_3)$.
The result is immediate from \fullref{main1}, \fullref{main4} 
and \fullref{main3}.

(4)\qua In the case $p=3$, $G=E_7$, 
the Weyl group $W(A)$ is the subgroup of $GL_4(\mathbb{F}_3)$ 
consisting of the following matrices:
\[
\left( \begin{array}{c|c}
\varepsilon_1 & m \\ \hline
0 & g_1 \end{array}
\right), 
\]
where $\varepsilon_1\in \mathbb{F}_3^{\times}=\{1, 2\}$, $m\in 
M_{1,3}(\mathbb{F}_3)$, $g_1\in SL_3(\mathbb{F}_3)$.
Firstly, we compute the ring of invariants of a subgroup $W_0$ 
of $W(A)$.
The subgroup $W_0$
is the subgroup of $W(A)$ consisting of the matrices
\[
\left( 
\begin{array}{c|c}
1 & m \\ \hline
0 & g_1
\end{array}
\right),
\]
where $m \in M_{1,3}(\mathbb{F}_3)$, $g_1 \in 
SL_{3}(\mathbb{F}_3)$.
By \fullref{main1}, \fullref{main4} and \fullref{main3}, we have
\[
(H^{*}BA)^{W_0}=\mathbb{F}_3[x_{26}, x_{36}, x_{48}, x_{54}] 
\otimes \mathbb{F}_3\{ 1, Q_I u_3, Q_Ju_4 \},
\]
where $I$ ranges over $A_{3}'$ and $J$ ranges over $A_3$.
Let 
\[
R= \mathbb{F}_3[x_{26}, x_{36}, x_{48}, x_{108}], 
\]
and let
\[
M=\mathbb{F}_3\{ x_{54}^{\delta}, x_{54}^{\delta}Q_I u_3, 
x_{54}^{\delta}Q_J u_4 \},
\]
where $x_{108}=x_{54}^2$,  $\delta\in \{0,1\}$, $I$ ranges over 
$A_{3}'$ and $J$ ranges over $A_3$.
Then, we have
\[
(H^{*}BA)^{W_0}=R\otimes M.
\]
Next, we calculate the ring of invariants $(H^{*}BA)^{W(A)}$ as 
a subspace of $(H^{*}BA)^{W_0}$.
Put
\[
\alpha=\left( \begin{array}{c|ccc}
2 & 0 & 0 & 0 \\ \hline
0 & 1 & 0 & 0 \\
0 & 0 & 1 & 0 \\
0 & 0 & 0 & 1
\end{array}
\right).
\]
Then, for $x\in R$, we have $\alpha x=x$ and we also have the 
following direct sum decomposition:
\[
M=M_1\oplus M_2,
\]
where $M_i =\{ x \in M | \alpha x=i x\}$ for $i=1$, $2$. 
In particular, we have
\[
M_1=\mathbb{F}_3\{ 1, Q_Iu_3, x_{54}Q_Ju_4 \}.
\]
Since the Weyl group $W(A)$ is generated by $W_0$ and $\alpha$, 
an element $x$ in $R\otimes M$
is $W(A)$--invariant if and only if $\alpha x=x$.
Hence, we have
\[
(H^{*}BA)^{W(A)}=R \otimes M_1.
\]

(5)\qua In the case $p=3$, $G=E_8$, $A=E_{E_8}^{5a}$, 
the Weyl group $W(A)$ is the subgroup of $GL_{5}(\mathbb{F}_3)$ 
consisting of the following matrices:
\[
\left(
\begin{array}{c|c|c}
\varepsilon_1 & m_0 & m_1 \\ \hline
0 & g_1 & 0 \\ \hline
0 & 0 & \varepsilon_2 
\end{array}
\right),
\]
where $\varepsilon_1$, $\varepsilon_2\in 
\mathbb{F}_3^{\times}=\{1,2\}$, $m_0\in M_{1,3}(\mathbb{F}_3)$, 
$m_1\in M_{1,1}(\mathbb{F}_3)$, $g_1\in SL_3(\mathbb{F}_3)$.
We consider the subgroup $W_0$ of $W(A)$ consisting of the 
following matrices:
\[
\left(
\begin{array}{c|c|c}
1 & m_0 & m_1 \\ \hline
0 & g_1 & 0 \\ \hline
0 & 0 & 1 
\end{array}
\right),
\]
where $g_1\in SL_3(\mathbb{F}_3)$, $m_0\in 
M_{1,3}(\mathbb{F}_3)$, $m_1 \in M_{1,1}(\mathbb{F}_3)$.
By \fullref{main4} and \fullref{main3}, we have
\[
(H^{*}BA)^{W_0}=\mathbb{F}_3[x_2, x_{26}, x_{36}, x_{48}, 
x_{162}]\otimes \mathbb{F}_3\{ 1, Q_I u_3, u_1, (Q_I u_3)u_1, 
Q_J u_5\},
\]
where $I$ ranges over $A'_3$ and $J$  ranges over $A_{4}$.
Let 
\[
R= \mathbb{F}_3[x_4, x_{26}, x_{36}, x_{48}, x_{324}], 
\]
and let
\[
M= \mathbb{F}_3\{ x_{2}^{\delta_1}x_{162}^{\delta_2}, 
x_{2}^{\delta_1}x_{162}^{\delta_2}Q_I u_3, 
x_{2}^{\delta_1}x_{162}^{\delta_2}u_1, 
x_{2}^{\delta_1}x_{162}^{\delta_2}(Q_I u_3)u_1, 
x_{2}^{\delta_1}x_{162}^{\delta_2}Q_J u_5\},
\]
where $\delta_1$, $\delta_2\in \{0,1\}$, $I$ ranges over $A'_3$ 
and $J$  ranges over $A_{4}$.
Consider matrices
\[
\alpha =\left(
\begin{array}{c|ccc|c}
2 & 0 & 0 & 0 & 0 \\ \hline
0 & 1 & 0 & 0 & 0 \\
0 & 0 & 1 & 0 & 0 \\
0 & 0 & 0 & 1 & 0 \\ \hline
0 & 0 & 0 & 0 & 1
\end{array}
\right), \qua
\beta =\left(
\begin{array}{c|ccc|c}
1 & 0 & 0 & 0 & 0 \\ \hline
0 & 1 & 0 & 0 & 0 \\
0 & 0 & 1 & 0 & 0 \\
0 & 0 & 0 & 1 & 0 \\ \hline
0 & 0 & 0 & 0 & 2
\end{array}
\right).
\]
Then, we have $\alpha x=x$ and $\beta x = x$ for $x\in R$ . 
Furthermore, it is also clear that we have the following direct 
sum decomposition:
\[
M=M_{1,1}\oplus M_{1,2} \oplus M_{2,1} \oplus M_{2,2},
\]
where 
\[
M_{i,j}=\{ x \in M | \alpha x = i x , \beta x = j x \}.
\]
In particular, we have
\[
M_{1,1}=\mathbb{F}_3\{ 1, x_2 u_1, Q_I u_3, x_2 (Q_I u_3)u_1, 
x_2x_{162} Q_Ju_5\}.
\]
Since $W(A)$ is generated by $W_0$ and $\alpha$, $\beta$ in the 
above, $x\in (H^{*}BA)^{W_0}$
is $W(A)$--invariant if and only if $\alpha x=\beta x =x$.
Hence, we have
\[
(H^{*}BA)^{W(A)}=R \otimes M_{1,1}.
\proved
\]
\end{proof}
 
\begin{rem}
Our computation of the ring of invariants of polynomial tensor 
exterior algebra in \fullref{ring}
is based on the computation of the ring of invariants of 
polynomial algebra 
and the assumption that the
ring of invariants of polynomial algebra is also  a polynomial 
algebra. 
In the case $A=E_{E_8}^{5b}$, however, 
the Weyl group does not satisfy the condition we assume in this 
section
and  the ring of invariants of polynomial algebra is no longer a 
polynomial algebra.
Hence, both \fullref{main4} and \fullref{main3} do not apply in 
this case.
\end{rem}

\section[O_n-1 and D_n-1]
{$\mathcal{O}_{n-1}(x_i)$ and $\mathcal{D}_{n-1}$} \label{section5}

In this section, we collect some facts, which we need in the 
proof of \fullref{main4} and \fullref{main3}.

For $i=1, \ldots, n$, 
the element $\mathcal{O}_{n-1}(x_i)$ in $\mathbb{F}_q[x_1, 
\ldots, x_{n}]$ is defined to be
\[
\mathcal{O}_{n-1}(x_i)=\prod_{x \in \mathbb{F}_q\{ x_2, \ldots, 
x_n\}} (x_i+x).
\]
We also define $\mathcal{O}_{n-2}(x_i)$ in $\mathbb{F}_q[x_1, 
\ldots, x_{n-1}]$ by
\[
\mathcal{O}_{n-2}(x_i)=\prod_{x \in \mathbb{F}_q\{x_2, \ldots, 
x_{n-1}\}} (x_i+x)
\]
for $n\geq 3$ and by
\[
\mathcal{O}_0(x_i)=x_i
\]
for $n=2$.

Using the same argument as in the  proof of \fullref{sum2}, 
we can easily obtain the following proposition.

\begin{prop} \label{sum5}
For $i=1, \ldots, n$, 
we may express $\mathcal{O}_{n-1}(x_i)$ and  
$\mathcal{O}_{n-2}(x_i)$ in terms of Dickson invariants 
as follows{\rm :}
\[
\begin{array}{ll}
\displaystyle \mathcal{O}_{n-1}(x_i)=\sum_{j=0}^{n-1} 
(-1)^{n-1-j}c_{n-1,j}(x_2, \ldots, x_n) {x_i}^{q^j} 
& \mbox{for $n\geq 2$},
\\
\displaystyle \mathcal{O}_{n-2}(x_i)=\sum_{j=0}^{n-2} 
(-1)^{n-2-j}c_{n-2,j}(x_2, \ldots, x_{n-1}) {x_i}^{q^j} 
&\mbox{for $n>2$}.
\end{array}
\]
\end{prop}

We need \fullref{prop52} and \fullref{prop53} below in the 
proof of \fullref{prop66}.

\begin{prop} \label{prop52}
In $\mathbb{F}_q[x_1, \ldots, x_{n-1}]$, we have the following 
equality{\rm :}
\[
e_{n-1}(x_1, \ldots, x_{n-1})=\mathcal{O}_{n-2}(x_1) 
e_{n-2}(x_2, \ldots, x_{n-1})
\]
for $n\geq 3$, and 
\[
e_1(x_1) = x_1
\]
for $n=2$.
\end{prop}

\begin{proof}
For $n=2$, the proposition is obvious.
For $n\geq 3$, by \fullref{divisible2}, \fullref{sum2} and 
\fullref{sum5}, we have
\[
\begin{array}[b]{rcl}
&&
\displaystyle e_{n-1}(x_1, \ldots, x_{n-1}) \\[1ex]
& = & Q_0\ldots Q_{n-2}dx_1\ldots dx_{n-1} \\
& = & \displaystyle \sum_{j=0}^{n-2} (-1)^{n-2-j} (Q_0\ldots 
\widehat{Q}_j \ldots Q_{n-2} dx_2 \ldots dx_{n-1}) x_1^{q^j}\\
& = & \displaystyle \sum_{j=0}^{n-2} (-1)^{n-2-j} e_{n-2}(x_2, 
\ldots, x_{n-1})c_{n-2,j}(x_2, \ldots, x_{n-1}) x_1^{q^j}\\
& = & \displaystyle e_{n-2}(x_2, \ldots, x_{n-1}) 
\mathcal{O}_{n-2}(x_1).
\end{array}
\proved
\]
\end{proof}

\begin{prop}\label{prop53}
The obvious projection
\[
\pi\co\mathbb{F}_q[x_1, \ldots, x_n]\longrightarrow 
\mathbb{F}_q[x_1, \ldots, x_{n-1}]
\]
maps $\mathcal{O}_{n-1}(x_1)$, $e_{n-1}(x_2, \ldots, x_n)$ to 
$\mathcal{O}_{n-2}(x_1)^q$, $0$, respectively. 
\end{prop}

\begin{proof}
Since $e_{n-1}(x_2, \ldots, x_{n})$ is divisible by $x_n$, 
we have 
\[
\pi(e_{n-1}(x_2, \ldots, x_n))=0
\]
as in the proof of \fullref{projection}.
For $n=2$, the equality 
\[
\pi(\mathcal{O}_1(x_1))=\mathcal{O}_0(x_1)^{q}=x_1^q
\]
is obvious. 
For $n\geq 3$,
by \fullref{projection}, 
we have 
\[
\pi(c_{n-1,j}(x_2, \ldots, x_n))=c_{n-2,j-1}(x_2, \ldots, 
x_{n-1})^q
\]
for $j=1, \ldots, n-1$. 
Hence, we have
\[
\begin{array}[b]{rcl}
\pi (\mathcal{O}_{n-1}(x_1)) & = & \displaystyle 
\pi \left( \sum_{j=0}^{n-1} (-1)^{n-1-j} c_{n-1, j}(x_2, \ldots, 
x_n) x_1^{q^{j}} \right) \\
& = &  \displaystyle \sum_{j=1}^{n-1} (-1)^{n-1-j} c_{n-2, 
j-1}(x_2, \ldots, x_{n-1})^q x_1^{q^j}  \\
& = & \displaystyle \mathcal{O}_{n-2}(x_1)^q.
\end{array} 
\proved
\]
\end{proof}

%: D_{n-1}

For $a$ in $P_n\otimes E_n$, let 
\[
\mathcal{D}_{n-1}(a)=\sum_{j=0}^{n-1} (-1)^{n-1-j}c_{n-1, 
j}(x_2, \ldots, x_n)Q_{j} a.
\]
Then, $\mathcal{D}_{n-1}$ induces a $P_n$--linear homomorphism
\[
\mathcal{D}_{n-1}\co P_n\otimes E_n^{r} \longrightarrow 
P_n\otimes E_{n}^{r-1},
\]
which extends naturally to 
\[
\mathcal{D}_{n-1}\co K_n \otimes E_{n}^{r}\to K_n\otimes 
E_{n}^{r-1}.
\]

\begin{prop} \label{d}
For $i=1, \ldots, n$, we have
\[
\mathcal{D}_{n-1}(dx_i)=\mathcal{O}_{n-1}(x_i).
\]
In particular, for $i=2, \ldots, n$, we have 
$\mathcal{D}_{n-1}(dx_i)=0$.
\end{prop}

\begin{proof}
By \fullref{sum5}, we have
\[
\begin{array}{rcl}
\mathcal{D}_{n-1}(dx_i)&=&\displaystyle \sum_{j=0}^{n-1} 
(-1)^{n-1-j}c_{n-1,j}(x_2, \ldots, x_n) Q_j dx_i \\
& = & \displaystyle \sum_{j=0}^{n-1} (-1)^{n-1-j} c_{n-1,j}(x_2, 
\ldots, x_n) x_i^{q^j} \\
& = & \displaystyle \mathcal{O}_{n-1}(x_i).
\end{array}
\]
On the other hand, by the definition of 
$\mathcal{O}_{n-1}(x_i)$, we have
\[
\mathcal{O}_{n-1}(x_i)=0
\]
for $i=2, \ldots, n$. Hence, we have $\mathcal{D}_{n-1}(dx_i)=0$ 
for $i=2, \ldots, n$.
\end{proof}

Let $g_1$ be an element in $GL_{n-1}(\mathbb{F}_q)$. 
We consider the following matrix
\[
\bar{g}_1=\left(\begin{array}{c|c} 
1 & 0 \\ \hline
0 & g_1
\end{array}
\right).
\]

We need \fullref{prop55} and \fullref{prop56} below in the 
proof of \fullref{prop68}.

\begin{prop}\label{prop55}
For $g_1$ in $GL_{n-1}(\mathbb{F}_q)$, there holds 
\[
\mathcal{D}_{n-1}(\bar{g
}_1a)=\bar{g}_1\mathcal{D}_{n-1}(a).
\]
\end{prop}

\begin{proof}
Suppose that 
\[
a = \sum_{i=1}^{n} a_i dx_i,
\]
where $a_i \in K_n$ for $i=1, \ldots, n$.
Since 
\[
\bar{g}_1 dx_1=dx_1,
\]
and since, for $i=2, \ldots, n$, 
\[
\bar{g}_1 dx_i
\]
is in $E_{n-1}^{1}$, we have
\[
\bar{g}_1 a = (\bar{g}_1a_1) dx_1 + a_2' dx_2+\cdots +a_n' dx_n
\]
for some $a_i'\in K_n$ for $i=2, \ldots, n$. 
Hence, by \fullref{d}, we have
\[
\mathcal{D}_{n-1}(\bar{g}_1 a)
=(\bar{g}_1a_1) \mathcal{O}_{n-1}(x_1).
\]
On the other hand, by \fullref{d}, we have
\[
\bar{g}_1 \mathcal{D}_{n-1}(a)= \bar{g}_1 (a_1 
\mathcal{O}_{n-1}(x_1))=(\bar{g}_1a_1)\mathcal{O}_{n-1}(x_1).
\proved
\]
\end{proof}

Let $\mathbb{P}^{n-2}$ be the projective space 
\[
(\mathbb{F}_q^{n-1}\backslash \{(0,\ldots,0)\})/\sim,
\]
where $\ell\sim \ell'$ if and only if there is $\alpha \in 
\mathbb{F}_q^{\times}$ such that $\ell=\alpha \ell'$.

\begin{prop} \label{prop57}\label{prop56}
If $a\in P_n$ is divisible by $\alpha_2x_2+\cdots+\alpha_nx_n$
for arbitrary $(\alpha_2, \ldots, \alpha_n) \allowbreak \in 
\mathbb{P}^{n-2}$, 
then $a$ is divisible by $e_{n-1}(x_2, \ldots, x_{n})$.
\end{prop}

\begin{proof}
Since the number of elements in $\mathbb{P}^{n-2}$ is equal to 
the homogeneous degree $1+q+\cdots+q^{n-2}$ 
of $e_{n-1}(x_2, \ldots, x_n)$,
it suffices to show that $e_{n-1}(x_2, \ldots, x_n)$ is a 
product of elements of the form 
$\alpha_2x_2+\cdots+\alpha_nx_n$, where $(\alpha_2, \ldots, 
\alpha_n)$ ranges over $\mathbb{P}^{n-2}$.
It is clear that $e_{n-1}(x_2, \ldots, x_n)$ is divisible by 
$x_n$. 
So, we have 
\[
e_{n-1}(x_2, \ldots, x_n)=b x_n
\]
for some $b\in P_n$. 
It is also clear that $e_{n-1}(x_2, \ldots, x_n)$ is invariant 
under the action of $SL_{n-1}(\mathbb{F}_q)$.
There is $g_1\in SL_{n-1}(\mathbb{F}_q)$ such that 
$\alpha_2x_2+\cdots +\alpha_nx_n=\bar{g}_1 x_n$. Hence, 
on the one hand, we have
\[
\bar{g}_1 e_{n-1}(x_2, \ldots, x_n)=(\bar{g}_1 b) 
(\alpha_2x_2+\cdots+\alpha_nx_n)
\]
and, on the other hand, we have
\[
\bar{g}_1 e_{n-1}(x_2, \ldots, x_n)=e_{n-1}(x_2, \ldots, x_n).
\]
Therefore, $e_{n-1}(x_2, \ldots, x_n)$ is divisible by arbitrary 
$\alpha_2x_2+\cdots+\alpha_nx_n$. This completes the proof.
\end{proof}

\section[Proof of Theorem 4.1]{Proof of \fullref{main4}} \label{section6}

In order to prove \fullref{main4}, 
we recall the strategy to compute rings of invariants given by 
Wilkerson in \cite[Section~3]{wilkerson}.
It can be stated in the following form.

\begin{thm}\label{wilkerson}
Suppose that  $G$ is a subgroup of $GL_{n}(\mathbb{F}_q)$ and 
$G$ 
acts on the polynomial algebra $\mathbb{F}_q[x_1, \ldots, x_n]$ 
in the obvious manner.
Let $f_1, \ldots, f_n$ be homogeneous $G$--invariant polynomials 
in $\mathbb{F}_q[x_1, \ldots, x_n]$.
Let $R$ be the subalgebra of $\mathbb{F}_q[x_1, \ldots, x_n]$ 
generated by $f_1, \ldots, f_n$.
Then, $R$ is a polynomial algebra $\mathbb{F}_q[f_1, \ldots, 
f_n]$ and the ring of invariants $\mathbb{F}_q[x_1,\ldots, 
x_n]^{G}$ is equal to 
the subalgebra $R$
if and only if
$\mathbb{F}_q[x_1,\ldots, x_n]$ is integral over $R$ and $\deg 
f_1\ldots \deg f_n=|G|$.
\end{thm}

In the statement of \fullref{wilkerson},
$\deg f$ is the homogeneous degree of $f$, that is, we define
the degree $\deg x_i$ of indeterminate $x_i$ to be $1$.
For the proof of this theorem, 
we refer the reader to
Smith's book \cite[Corollaries 2.3.2 and 5.5.4, and 
Proposition 5.5.5]{smith} and Wilkerson's paper 
\cite[Section 3]{wilkerson}.

\begin{proof}[Proof of \fullref{main4}]

As we mentioned, in order to prove \fullref{main4}, it suffices 
to show the following:
\begin{itemize}
\item[(1)]  homogeneous polynomials $
\displaystyle \mathcal{O}_{n-1}(x_1), f_2, \ldots, f_n$
are $G$--invariant{\rm ;}
\item[(2)] indeterminates $\displaystyle x_1, \ldots, x_n$
are integral over $R${\rm ;}
\item[(3)]
the product of homogeneous degrees of $\displaystyle 
\mathcal{O}_{n-1}(x_1), f_2, \ldots, f_n$
is equal to the order of $G$, that is, 
\[
\displaystyle \deg \mathcal{O}_{n-1}(x_1)\deg f_2\ldots \deg 
f_n=|G|.
\]
\end{itemize}

By definition, 
$f_2, \ldots, f_n$ are $G_1$--invariant, and so they are also 
$G$--invariant.
It follows from \fullref{wilkerson} that
$x_2,\ldots, x_n$ are integral over 
$R_1=\mathbb{F}_q[f_2,\ldots, f_n]$, and so they are integral 
over $R$.
It is also immediate from \fullref{wilkerson} that $\deg f_2\ldots\deg f_n=|G_1|$.
It is clear from the definition of $\mathcal{O}_{n-1}(x_1)$  that
$\deg \mathcal{O}_{n-1}=2^{n-1}$.
Hence, we have  
\[
\deg \mathcal{O}_{n-1}(x_1)\deg f_2\ldots \deg 
f_n=2^{n-1}|G_1|=|G|.
\]

So, it remains to show the following:
\begin{itemize}
\item[(1)] 
$\mathcal{O}_{n-1}(x_1)$ is $G$--invariant and
\item[(2)] $x_1$ is integral over $R$.
\end{itemize}

First, we deal with (1). By the definition of 
$\mathcal{O}_{n-1}(x_1)$, we have that
\[
g \mathcal{O}_{n-1}(x_1)
\]
is a product of
\[
g  (x_1+x)=x_1+\sum_{j=2}^{n} a_{1, j}(g^{-1})x_j +gx,
\]
where $x$ ranges over $\mathbb{F}_q\{x_2, \ldots, x_n\}$. 
As $x$ ranges over 
$\mathbb{F}_q\{x_2,\ldots, x_n\}$, the sum
\[
\sum_{j=2}^{n} a_{1, j}(g^{-1})x_j +gx
\]
 also ranges over 
$
\mathbb{F}_q\{x_2,\ldots, x_n\}
$.
Hence, we have 
\[
g \mathcal{O}_{n-1}(x_1)=\mathcal{O}_{n-1}(x_1).
\]

Next, we deal with (2). 
By \fullref{sum5}, we have
\[
\mathcal{O}_{n-1}(X)=X^{q^{n-1}}+\sum_{j=0}^{n-2}(-1)^{n-1-j} 
c_{n-1,j}(x_2, \ldots, x_n) X^{q^{j}}.
\]
Since Dickson invariants $c_{n-1,j}(x_2, \ldots, x_n)$ are in 
$R_1=\mathbb{F}_2[x_2, \ldots, x_n]^{G_1}$, 
the polynomial 
\[
\varphi(X)=\mathcal{O}_{n-1}(X) -\mathcal{O}_{n-1}(x_1)
\]
is a monic polynomial in $R[X]$.
It is clear that \[
\varphi(x_1)=0.
\]
Hence, the indeterminate $x_1$ is integral over $R$.
This completes the proof.
\end{proof}

\section[Proof of Theorem 4.2]{Proof of \fullref{main3}} \label{section7}

Let $G_0$ be the subgroup of $G$ consisting of the following 
matrices:
\[
\left( \begin{array}{c|c}
1 &m \\ \hline
0 &   1_{n-1}  \end{array}\right),
\]
where $m\in M_{1,n-1}(\mathbb{F}_q)$, $1_{n-1}$ is the identity 
matrix in $GL_{n-1}(\mathbb{F}_q)$. 
Let $B_{n}$ be the set of subsets of 
\[
\{ 2, \ldots, n\}.
\]
Let $B_{n,r}$ be the subset of $B_n$
such that
$J\in B_{n,r}$ if and only if 
\[
J=\{j_1, \ldots, j_r\} \qua \mbox{ and } \qua 1< 
j_1<\cdots<j_{r}\leq n.  \]
We write $dx_{J}$ for \[
dx_{j_1}\ldots dx_{j_r}
\] and we define $dx_{\emptyset}$ to be $1$.

The following proposition is nothing but the particular case of 
\fullref{main4} and \fullref{main3}.

\begin{prop} \label{prop0}
The ring of invariants $P_{n}^{G_0}$ is given as follows{\rm :}
\[
P_{n}^{G_0}=\mathbb{F}_q[\mathcal{O}_{n-1}(x_1), x_2, \ldots, 
x_n].
\]
The ring of invariants $(P_{n}\otimes E_{n})^{G_0}$ is a free 
$P_{n}^{G_0}$--module with the basis
\[
\{ Q_{I}dx_1\ldots dx_n, dx_{J}\},
\]
where $I$ ranges over $A_{n-1}$ and  $J$ ranges over $B_{n}$.
\end{prop}

Now, we consider a $K_n$--basis for $K_n \otimes E_n$.

\begin{prop} \label{prop1}
The elements
\[
Q_{0}\ldots Q_{n-2}dx_1\ldots dx_n, dx_2, \ldots, dx_n
\] form a $K_n$--basis for $K_n\otimes E_n^{1}$.
\end{prop}

\begin{proof}
For dimensional reasons, it suffices to show that the above 
elements are linearly independent in $K_n\otimes E_n^1$.
Suppose that
\[
a_1 Q_0\ldots Q_{n-2}dx_1\ldots dx_n +a_2dx_2+\cdots +a_n dx_n=0
\]
in $K_n\otimes E_n^1$, where $a_1, \ldots, a_n$ are in $K_n$.
Then, since 
\[
Q_0 \ldots Q_{n-2}dx_1\ldots dx_n=\sum_{i=1}^{n} (-1)^{n-i} 
e_{n-1}(x_1, \ldots, \widehat{x}_i, \ldots, x_n)dx_i,
\]
we have
\begin{eqnarray*}
(-1)^{n-1}a_1 e_{n-1}(\widehat{x}_1, \ldots, x_n)& = & 0, \\
a_2+(-1)^{n-2} a_1e_{n-1}(x_1,\widehat{x}_2, \ldots, x_n) & = & 0, \\
&\vdots&\\
a_n +(-1)^0 a_1 e_{n-1}(x_1, \ldots, x_{n-1}, \widehat{x}_n)& = & 0. 
\end{eqnarray*}
Thus, solving this linear system, we obtain $a_1=0, \ldots, 
a_n=0$.
\end{proof}

\begin{prop} \label{prop2}
The elements
\[
Q_{I}dx_1\ldots dx_n, dx_{J}
\] form a $K_n$--basis for $K_n\otimes E_n^{r}$,
where $I\in A_{n-1, n-r}$ and $J\in B_{n,r}$.
\end{prop}

\begin{proof}
Again, for dimensional reasons, it suffices to show that the 
above elements are linearly independent
in $K_n\otimes E_n^r$.
Suppose that 
\[
\sum_{I\in A_{n-1, n-r}} a_{I} Q_{I} dx_1\ldots dx_n + 
\sum_{J\in B_{n,r}} b_{J}dx_{J}=0,
\]
where $a_{I}$, $b_{J}$ are in $K_{n}$. 
The linear independence of the terms $dx_{J}$ is clear.  Hence, it 
remains to show that $a_{I}=0$ for each $I$.

Fix $I \in A_{n-1,n-r}$ and let $K=S_{n-1}\backslash I$.
Then, applying $Q_{K}$ to the both sides of the above equality, 
we have
\[
\mathrm{sign} (K,I) a_{I}Q_{0}\ldots Q_{n-2} dx_1\ldots 
dx_n+\alpha=0
\]
in $K_{n} \otimes E_{n}^{1}$, where $\alpha$ is a linear 
combination of $dx_2, \ldots, dx_n$ over $K_n$.
Hence, by \fullref{prop1}, we have $a_{I}=0$ for each $I \in 
A_{n-1, n-r}$.
\end{proof}

Suppose that $a$ is in $P_n \otimes E_n^1$.
Then, on the one hand, 
we may express $a$ as follows:
\[
a = \varphi_1 dx_1 +\cdots +\varphi_n dx_n,
\]
where $\varphi_1, \ldots, \varphi_n$ are in $P_n$.
On the other hand, by \fullref{prop1}, 
we may express $a$ as follows:
\[
a = a_1 Q_0\ldots Q_{n-2}dx_1\ldots dx_n + a_2dx_2+\cdots 
+a_ndx_n,
\]
where $a_1, \ldots, a_n$ are in $K_n$. Observe that terms $a$ 
and $\varphi$ are unique in the above expressions.

We need to show that $a_1, \ldots, a_n$ are in $P_n$ if $a$ is 
$G_0$--invariant.

\begin{prop} \label{propPoly}
There are polynomials $a_i'$ over $\mathbb{F}_q$ in $n$ 
variables such that
 \[
 a_ie_{n-1}(x_2, \ldots, x_n)=a_{i}'(x_1, \ldots, x_n)
 \]
 for $i=1, \ldots, n$.
 \end{prop}
 
\begin{proof}
For $i=1$,   
we apply $\mathcal{D}_{n-1}$ to $a$.
Then, we have
\[
\mathcal{D}_{n-1}(a)=\varphi_1 \mathcal{O}_{n-1}(x_1).
\]
On the other hand, we have 
\[
\begin{array}{rcl}
\mathcal{D}_{n-1}(a)&=&(-1)^{n-1}a_1Q_0\ldots Q_{n-1}dx_1\ldots 
dx_n\\
& =&(-1)^{n-1}a_{1}\mathcal{O}_{n-1}(x_1)e_{n-1}(x_2, \ldots, 
x_n).
\end{array}
\]
Hence, we obtain
\[
a_{1}=(e_{n-1}(x_2, \ldots, x_n))^{-1} \varphi_1
\]
and 
\[
a_{1}'(x_1, \ldots, x_n)=\varphi_1.
\]
For $i=2, \ldots, n$, 
applying $Q_0, \ldots, Q_{n-2}$ to $a$, we have a linear system
\begin{eqnarray*}
Q_0 a &  = & a_2 x_2 + \cdots + a_n x_n, \\
Q_1 a & = & a_2 x_2^{q}+\cdots + a_n x_n^{q}, \\
 & \vdots & \\
Q_{n-2} a & = & a_2 x_{2}^{q^{n-2}} +\cdots + a_n x_{n}^{q^{n-2}}.
\end{eqnarray*}
Writing this linear system in terms of matrix, 
we have
\[
\left( \begin{array}{c} Q_0 a \\ Q_1a \\ \vdots \\ Q_{n-2}a 
\end{array} \right) 
=A \left(
\begin{array}{c}
a_2 \\ a_3 \\ \vdots \\ a_n  \end{array}
\right),
\]
where
\[
A =\left( \begin{array}{ccc}
x_2 & \ldots & x_n  \\
x_2^q & \ldots & x_n^q \\
\vdots & \ddots & \vdots  \\
x_2^{q^{n-2}} & \ldots & x_{n}^{q^{n-2}} \end{array}
\right).
\]
It is clear that $\det A=e_{n-1}(x_2, \ldots, x_n)\not=0$.
It is also clear  that each entry of $A$ is in $P_{n-1}$.
Therefore,  for some $\varphi_{i,j}$ in $P_{n-1}$, we have
\[
a_i = e_{n-1}(x_2, \ldots, x_n)^{-1}\left(\sum_{j=1}^{n} 
\varphi_{i,j} Q_{j}a \right).
\]
Since $Q_{j}a$ is in $P_{n}$, by letting 
\[ a_{i}'(x_1, \ldots, x_n)=\sum_{j=1}^{n} \varphi_{i,j} Q_{j}a, 
\]
we obtain the required results.
\end{proof}

\begin{prop}\label{prop66}
Suppose that $a$ is ${G_0}$--invariant. Then
$a_i'(x_1, \ldots, x_n)$ in \fullref{propPoly} are also 
${G_0}$--invariant
for $i=1, \ldots, n$.
\end{prop}

\begin{proof}
For $g \in G_0$, we have
\[
g (dx_i)=dx_i
\]
for $i=2, \ldots, n$
and, since $G_0\subset SL_n(\mathbb{F}_q)$, we have
\[
g (dx_1\ldots dx_n)= dx_1\ldots dx_n.
\]
Since the action of Milnor operations $Q_j$ commutes with the 
action of 
the general linear group $GL_n(\mathbb{F}_q)$, we have
\[
g(Q_0\ldots Q_{n-2}dx_1\ldots dx_n)=Q_0\ldots Q_{n-2}dx_1\ldots 
dx_n.
\]
Hence, we have
\[
g a = (g a_1) Q_0\ldots Q_{n-2}dx_1\ldots dx_n +(ga_2) 
dx_2+\cdots +(ga_n) dx_n.
\]
Thus, if $ga=a$, we have $ga_1=a_1, \ldots, ga_n=a_n$.
It is also clear that 
$e_{n-1}(x_2, \ldots, x_n)$ is ${G_0}$--invariant.
Hence, $a_{i}'(x_1, \ldots, x_n)$ are also ${G_0}$--invariant for 
$i=1, \ldots, n$.
\end{proof}

\begin{prop}
Suppose that $a$ is ${G_0}$--invariant. Then, 
$a_{1}'(x_1, \ldots. x_n)$ is divisible by $x_n$.
\end{prop}

\begin{proof}
Let us consider the coefficient $\varphi_n$ of $dx_n$;
we have 
\[
\varphi_n e_{n-1}(x_2, \ldots, x_n) = {a}'_{n}(x_1, \ldots, x_n) 
+ a_1'(x_1, \ldots, x_n) e_{n-1}(x_1, \ldots, x_{n-1}) .
\]
Since $a_{1}'(x_1, \ldots, x_n)$ and $a_{n}'(x_1, \ldots, x_n)$ 
are 
$G_0$--invariant, 
there are polynomials $a_{1}''$, $a_{n}''$ over $\mathbb{F}_q$ 
in $n$ variables such that
\[
\begin{array}{l}
a_{1}'(x_1, \ldots, x_n)=a_{1}''(\mathcal{O}_{n-1}(x_1), x_2, 
\ldots, x_n), \\
a_{n}'(x_1, \ldots, x_n)=a_{n}''(\mathcal{O}_{n-1}(x_1), x_2, 
\ldots, x_n).
\end{array}
\]
Since $\mathcal{O}_{n-1}(x_1), x_2, \ldots, x_n$ are 
algebraically independent, 
it suffices to show that \[
a_{1}''(y_1, x_{2}, \ldots, x_{n-1}, 0)=0
\]
for algebraically independent $y_1, x_2, \ldots, x_{n-1}$.

Substituting $x_n=0$, we have the obvious projection 
\[
\pi\co\mathbb{F}_q[x_1, \ldots, x_n]\to \mathbb{F}_q[x_1, 
\ldots, x_{n-1}].
\]
It is clear from \fullref{prop52} and \fullref{prop53} that
\[
\begin{array}{rcll}
\pi(e_{n-1}(x_1, \ldots, x_{n-1})) &= &e_{n-1}(x_1, \ldots, 
x_{n-1})\\
& = &\mathcal{O}_{n-2}(x_1)e_{n-2}(x_2, \ldots, x_{n-1}) & 
\mbox{for $n\geq 3$}, \\
\pi(e_{1}(x_1)) &= & \mathcal{O}_{0}(x_1) & \mbox{for $n=2$}, \\
\pi(e_{n-1}(x_2, \ldots, x_{n})) & = & 0, \\
\pi(\mathcal{O}_{n-1}(x_1)) & = & \mathcal{O}_{n-2}(x_1)^{q}.
\end{array}
\]

Hence, for $n\geq 3$, we have
\[
0={a}_{n}''(y^q, x_2, \ldots, x_{n-1}, 0) +y{a}_1''(y^q, x_2, 
\ldots, x_{n-1},0)e_{n-2}(x_2, \ldots, x_{n-1}),
\]
where $y=\mathcal{O}_{n-2}(x_1)$ and $y, x_2, \ldots, x_{n-1}$ 
are algebraically independent.
Applying the partial derivative $\partial/\partial y$, we have
\[
{a}_1''(y^q, x_2, \ldots, x_{n-1},0)e_{n-2}(x_2, \ldots, 
x_{n-1})=0.
\]
Hence, we have
\[
a_{1}''(y^q, x_2, \ldots, x_{n-1}, 0)=0.
\]
Since $y^q, x_2,  \ldots, x_{n-1}$ are algebraically 
independent, 
we have the required result.

For $n=2$, we have
\[
0=a_{2}''(y^q, 0)+ya_{1}''(y^q,0).
\]
Applying the partial derivative $\partial/\partial y$, we have
\[
{a}_1''(y^q,0)=0
\]
and the required result.
\end{proof}

\begin{lem}\label{inv}
Suppose that $g_1\in GL_{n-1}(\mathbb{F}_q)$ and that $a$ is 
$G_0$--invariant.
 Then, $\bar{g}_1a$ is also $G_0$--invariant.
\end{lem}

\begin{proof}
Suppose that for each $g$, there is a $g'\in G_0$ such that 
$g\bar{g}_1=\bar{g}_1g'$. If it is true, then
for any $G_0$--invariant $a$, we have
\[
g\bar{g}_1a=\bar{g}_1g'a=g_1a.
\]
Hence, $g_1$ induces a homomorphism from $P_{n}^{G_0}$ to 
$P_{n}^{G_0}$. So,
it suffices to show that for each $g$ in $G_0$, there is a 
$g'\in G_{0}$ such that
$gg_1=g_1 g'$, which is
immediate from the following equality:
\[
\left( \begin{array}{c|c}
1 &   m  \\
\hline
0 & 1_{n-1}
\end{array}
\right)
\left( \begin{array}{c|ccc}
1 & 0 \\
\hline
0 &  g_1
\end{array}
\right)=
\left( \begin{array}{c|c}
1 & 0 \\
\hline
0 &  g_1 \end{array}
\right)
\left( \begin{array}{c|c}
1 &   mg_1  \\
\hline
0 & 1_{n-1}
\end{array}
\right),
\]
where $m\in M_{1, n-1}(\mathbb{F}_q)$ and $1_{n-1}$ stands for 
the identity matrix in $GL_{n-1}(\mathbb{F}_q)$.
\end{proof}

\begin{prop} \label{divisible6}\label{prop68}
Suppose that $a$ is $G_0$--invariant.
Then, 
$a_{1}, \ldots, a_n$ are in $P_{n}$.
\end{prop}

\begin{proof}
Firstly, we verify that $a_1$ is in $P_n$.
To this end, we prove that  the element $a'_1(x_1,\ldots, x_n)$ 
is divisible by $e_{n-1}(x_2, \ldots, x_n)$.
Let $\ell=\alpha_2x_2+\cdots+\alpha_nx_n$, where $\alpha_2, 
\ldots, \alpha_n\in \mathbb{F}_q$
 and $\ell\not =0$. 
 By \fullref{prop57}, it suffices to show that $a'_1(x_1, 
\ldots, x_n)$ is divisible by $\ell$.
There is $g_1$ in $GL_{n-1}(\mathbb{F}_q)$ such that 
$\bar{g}_1(x_n)=\ell$.
Since, by \fullref{inv},  $\bar{g}_1^{-1}a$ is also in 
$(P_n\otimes E_{n}^{1})^{G_0}$, there is an element $f$ in 
${P_{n}}$ such that
\[
\mathcal{D}_{n-1}(\bar{g}_1^{-1}a)=f x_n\mathcal{O}_{n-1}(x_1). 
\]
Here we have
\[
\mathcal{O}_{n-1}(x_1)a'_{1}(x_1,\ldots,x_n)=\mathcal{D}_{n-1}(a)
=\bar{g}_1 \mathcal{D}_{n-1}(\bar{g}_1^{-1}a)
=\mathcal{O}_{n-1}(x_1)\bar{g}_1(f) \ell. 
\]
So we have 
\[
a'_{1}(x_1,\ldots,x_n)=(\bar{g}_1f)\ell.
\]

Secondly, we verify that $a_i$ are in $P_n$ for $i=2, \ldots, 
n$, which follows
from the fact that
\[
\varphi_{i} = a_{i}+(-1)^{n-i} a_{1}e_{n-1}(x_1, \ldots, 
\widehat{x}_i, \ldots, x_n). 
\proved
\]
\end{proof}

\begin{proof}[Proof of \fullref{prop0}]

Suppose that
$a$ is in $P_{n}\otimes E_{n}^{r}$ and $G_0$--invariant.
By \fullref{prop2}, there are $a_{I}$, $b_J \in K_{n}$
such that
\[
a=\sum_{I\in A_{n-1, n-r}}a_{I}Q_{I}dx_1\ldots dx_n+\sum_{J\in 
B_{n,r}} b_{J}dx_J.
\]
It suffices to show that $a_I$, $b_J$ are in $P_n$.

Firstly, we verify that $a_{I}$ is in $P_{n}$.
Choose $I$ and let $K=S_{n-1}\backslash I$.
Then, we have
\[
Q_{K} a = \mathrm{sign}(K,I) a_I Q_{0}\ldots Q_{n-2} dx_1\ldots 
dx_n+\sum_{J\in B_{n,r}} b_{J} Q_{K} dx_J.
\]
By \fullref{divisible6}, $\mathrm{sign}(K,I) a_{I} $ is in 
$P_{n}$ and, by definition, $
\mathrm{sign}(K,I)\not=0$, hence $a_{I}$ is also in $P_{n}$. 

Secondly, we prove that $b_J$ is in $P_{n}$.
Put 
\[
a' = a-\sum_{I\in A_{n-1,n-r}} a_{I}Q_Idx_1\ldots 
dx_n=\sum_{J\in B_{n,r}} b_{J} dx_{J}.
\]
It is clear that $a'$ is also in $P_n \otimes E_n^{r}$. 
Hence, $b_J$ is in $P_n$. This completes the proof.
\end{proof}

%: Proof of Theorem

Now, we complete the proof of  \fullref{main3}.

\begin{proof}[Proof of \fullref{main3}]

Suppose that $a$ is an element in $P_n \otimes E_n$ and that $a$ 
is also $G$--invariant.
It suffices to show that $a$ is a linear combination of $\{ v_i, 
Q_I dx_1\ldots, dx_n\} $ over $P_n^{G}$.
It is clear that $a$ is also ${G_0}$--invariant. Hence, 
by \fullref{prop0}, 
there are $a_I$, $b_J$ in 
$P_{n}^{G_0}=\mathbb{F}_{q}[\mathcal{O}_{n-1}(x_1),x_2\ldots, 
x_n]$ such that
\[
a=\sum_{I} a_I Q_I dx_1\ldots dx_n+\sum_{J} b_{J} dx_J.
\]
Thus, we have
\[
a = \sum_{I} \sum_{k\geq 0} a_{I, k}\mathcal{O}_{n-1}(x_1)^k Q_I 
dx_1\ldots dx_n
+ \sum_{J} \sum_{k\geq 0} b_{J,k} \mathcal{O}_{n-1}(x_1)^{k} 
dx_J,
\]
where $a_{I,k}$, $b_{J,k}$ are in $P_{n-1}=\mathbb{F}_q[x_2, 
\ldots, x_n]$.
Since, by \fullref{main4}, $g\in G$ acts trivially on 
$\mathcal{O}_{n-1}(x_1)$, 
and since $g \in G\subset SL_{n}(\mathbb{F}_q)$ acts trivially 
on $Q_I dx_1\ldots dx_n$, we have 
\[
ga = \sum_{k\geq 0} \sum_{I} (ga_{I,k}) 
\mathcal{O}_{n-1}(x_1)^{k} Q_{I}dx_1\ldots dx_n
+ \sum_{k\geq 0}\sum_{J}  (g b_{J,k}) \mathcal{O}_{n-1}(x_1)^{k} 
(gdx_J). \]
It is clear that $gdx_J$ is in $E_{n-1}$.
As a $P_{n-1}$--module, $(P_{n}\otimes E_{n}^{r})^{G_0}$ is a 
free 
$P_{n-1}$--module with the 
basis 
\[
\{ \mathcal{O}_{n-1}(x_1)^k dx_J, 
\mathcal{O}_{n-1}(x_1)^{k}Q_{I}dx_1\ldots dx_n \}.
\]
Hence, we have
\[
g(a_{I,k})=a_{I,k}.
\]
Thus, $a_{I,k}$ is in ${P_{n-1}}^{G_1}$ and so $a_I$ is in 
$P_n^{G}$.
Put \[
a'=a-\sum_{I} a_{I} Q_{I}dx_1\ldots dx_n.
\]
Then, $a'$ is also in the ring of invariants $(P_n\otimes 
E_n)^G$, and we have
\[
a'=\sum_{k\geq 0} \left(\sum_{J} b_{J, k}dx_J 
\right)\mathcal{O}_{n-1}(x_1)^k.
\]
Hence, $\displaystyle \sum_{J}b_{J,k} dx_{J}$ is in the ring of 
invariants $(P_{n-1}\otimes E_{n-1})^{G_1}$.
By the assumption on the ring of invariants $(P_{n-1}\otimes 
E_{n-1})^{G_1}$,
there are polynomials $b_{1,k}, \ldots, b_{2^{n-1},k}$ in 
$P_{n-1}^{G_1}$  such that
\[
\sum_{J}b_{J,k} dx_{J}=\sum_{i=1}^{2^{n-1}} b_{i,k}v_{i}.
\]
Thus, writing $b_{i}$ for 
$\displaystyle \sum_{k\geq 0} b_{i,k} \mathcal{O}_{n-1}(x_1)^k$, 
we have
\[
a=\sum_{i=1}^{2^{n-1}} b_{i}v_i+
\sum_{I} a_{I} Q_{I}dx_1\ldots dx_n,
\]
where $b_i$, $a_{I}$ are in $P_{n}^{G}$.
This completes the proof.
\end{proof}

\appendix
\section[Appendix A]{Appendix} \setobjecttype{App}\label{sectionapp}

In \cite{mui}, 
M\`ui used the determinant of the $k\times k$ matrix 
$\bigl(\smash{x_j^{q^{i_\ell}}}\bigr)$ whose $(\ell,j)$ entry is~$\smash{x_j^{q^{i_\ell}}}$
to describe the Dickson invariant $[i_1, \ldots, i_k]$.
Using these Dickson invariants, he defined the M\`ui invariant 
$[r:i_1,\ldots, i_{n-r}]$. 
In this appendix, we verify in \fullref{amain} that the M\`ui 
invariant
\[
Q_{i_1}\ldots Q_{i_{n-r}}dx_1\ldots dx_n
\] in this paper
is indeed equal to the M\`ui invariant $[r: i_1,\ldots, i_{n-r}]$, 
up to sign.

Firstly, we recall the definitions. See \cite[Section 2]{mui} for 
the definition of $[i_1,\ldots, i_k]$ and
and \cite[Defintion 4.3]{mui} for the definition of 
$[r:i_1,\ldots, i_{n-r}]$.

\begin{defn}
The Dickson invariant
$[i_1, \ldots, i_k](x_1,\ldots, x_k)\in \mathbb{F}_q[x_1,\ldots, 
x_k]$ is 
defined by
\[
\sum_{\sigma}\mathrm{sgn} (\sigma) x_{\sigma(1)}^{q^{i_1}}\ldots 
x_{\sigma(k)}^{q^{i_k}}=\det (x_{j}^{q^{i_\ell}}),
\]
where $\sigma$ ranges over the set of permutations of 
$\{1,\ldots, k\}$ and $\mathrm{sgn}(\sigma)$ is the sign of the 
permutation $\sigma$.
\end{defn}

\begin{defn}
The M\`ui invariant 
$[r:i_1, \ldots, i_{n-r}]\in P_n\otimes E_n^{r}$ is defined by
\[
[r:i_1, \ldots, i_{n-r}]=\sum_{J} \mathrm{sgn}(\sigma_J) 
dx_{j_1}\ldots dx_{j_{r}} [i_1,\ldots, 
i_{n-r}](x_{j_{r+1}},\ldots, x_{j_n}),
\]
where 
\[
\sigma_J=\left( \begin{array}{ccc}
1, & \ldots, & n \\
j_1, & \ldots, & j_n
\end{array}
\right)
\]
ranges over the set of permutations of $\{1,\ldots, n\}$
such that  $j_1<\cdots<j_{r}$ and $j_{r+1}<\cdots<j_n$.
The above $\sigma_J$ corresponds to the subset $J=\{j_1, \ldots, 
j_r\}$ of order $r$  of $\{1, \ldots, n\}$.
\end{defn}

Secondly, we prove the following proposition.

\begin{prop} \label{ad}
There holds
\[
\begin{array}{rcl}
[i_1,\ldots, i_k](x_1,\ldots, x_k)&=&\displaystyle Q_{i_k}\ldots 
Q_{i_1}dx_1\ldots dx_k\\
&=&\displaystyle (-1)^{k(k-1)/2} Q_{i_1}\ldots Q_{i_k}dx_1\ldots 
dx_k.
\end{array}
\]
\end{prop}

\begin{proof}
We prove the first equality in this proposition by induction on 
$k$. Indeed, in the case $k=1$, the proposition holds.
Suppose that $k\geq 2$ and that there holds the equality
\[
[i_2, \ldots, i_k](x_1, \ldots, x_{k-1})=Q_{i_k}\ldots Q_{i_2} 
dx_1\ldots dx_{k-1}.
\]
Using the cofactor expansion (or the Laplace development) of the 
$k\times k$ matrix $(x_{j}^{q^{i_\ell}})$ along the first row 
\[
(x_1^{q^{i_1}}, x_{2}^{q^{i_1}}, \ldots, x_{k}^{q^{i_1}}),
\]
we have
\[
\begin{array}[b]{rcl}
[i_1, \ldots, i_k](x_1, \ldots, x_k)&=&\displaystyle 
\sum_{s=1}^{k} (-1)^{s+1} [i_2, \ldots, i_k](x_1, \ldots, 
\widehat{x}_{s}, \ldots, x_k)x_s^{q^{i_1}} \\
&=&\displaystyle  \sum_{s=1}^k (-1)^{s+1} (Q_{i_k}\ldots 
Q_{i_2}dx_1\ldots \widehat{dx}_{s}\ldots dx_k) x_{s}^{q^{i_1}} \\
&=&\displaystyle  Q_{i_k}\ldots Q_{i_2} \left( 
\sum_{s=1}^{k} (-1)^{s+1} x_{s}^{q^{i_1}}  dx_1\ldots 
\widehat{dx}_{s} \ldots dx_k \right) \\
&=&\displaystyle  Q_{i_k}\ldots Q_{i_1} dx_1\ldots dx_k. 
\end{array}
\]
So, the first equality holds. The second equality is immediate 
from the fact that
\[
Q_{i_k}\ldots Q_{i_1}=(-1)^{k(k-1)/2} Q_{i_1}\ldots Q_{i_k}. 
\proved
\]
\end{proof}

Finally, we state and prove the following proposition.

\begin{prop}\label{amain}
There holds
\[
\begin{array}{rcl}
[r:i_1,\ldots, i_{n-r}]&=&\displaystyle (-1)^{(n-r)r} 
Q_{i_{n-r}}\ldots Q_{i_1}dx_1\ldots dx_n\\
&=&\displaystyle (-1)^{(n-r)r+(n-r)(n-r-1)/2}Q_{i_1}\ldots 
Q_{i_{n-r}}dx_1\ldots dx_n.
\end{array}
\]
\end{prop}

\begin{proof}
As in the definition of $[r:i_1,\ldots, i_{n-r}]$, let 
$\sigma_J$ be a permutation of $\{1, \ldots, n\}$
with
$\sigma_{J}(1)<\cdots<\sigma_{J}(r)$, $\sigma_{J}(r+1)<\cdots 
<\sigma_{J}(n)$
and we denote by $j_k$ the value $\sigma_{J}(k)$ of $\sigma_J$ 
at $k$.
Let $I(J)$ be the ideal of $P_n \otimes E_n$ 
generated by $dx_{j_{r+1}}, \ldots, dx_{j_n}$.
Let \[
p_J\co
P_n\otimes E_n\to  P_n \otimes E_n/I(J)
\]
be the projection.
It is clear that 
\[
P_n \otimes E_n^r/((P_n \otimes E_n^r)\cap I(J))=P_n
\]
and 
for $f$ in $P_n\otimes E_n^r$, we have
\[
f=\sum_{J} p_{J}(f) dx_{j_1}\ldots dx_{j_r}.
\]
So, by \fullref{ad}, in order to prove the proposition, it 
suffices to show that 
\[
p_J(Q_{i_{n-r}}\ldots Q_{i_1} dx_1\ldots dx_n)
=
(-1)^{(n-r)r} \mathrm{sgn}(\sigma_J) Q_{i_{n-r}}\ldots Q_{i_1} 
dx_{j_{r+1}}\ldots dx_{j_n}.
\]
Suppose that 
\[
\psi(Q_{i_{n-r}}\ldots Q_{i_1})=1\otimes Q_{i_{n-r}}\ldots 
Q_{i_1}+ \sum a \otimes a',
\]
where $\psi$ is the coproduct of the Steenrod algebra.
We may choose $a'$, so that $a'=Q_{e_1}\ldots Q_{e_\ell}$ and  
$\ell<n-r$.
Thus, $a'(dx_{j_{r+1}}\ldots dx_{j_n})$ belongs to $I(J)$.

Then, there holds
\[
\begin{array}{rcl}
&&Q_{i_{n-r}}\ldots Q_{i_1} dx_1\ldots dx_n \\
&=& \displaystyle \mathrm{sgn} (\sigma_J) Q_{i_{n-r}}\ldots 
Q_{i_1} dx_{j_1}\ldots dx_{j_n} \\
&=& \displaystyle (-1)^{(n-r)r} \mathrm{sgn}( \sigma_J )
dx_{j_1}\ldots dx_{j_{r}} Q_{i_{n-r}}\ldots 
Q_{i_1}dx_{j_{r+1}}\ldots dx_{j_n}\\
&&\displaystyle +\sum
(-1)^{r\deg a'} \mathrm{sgn} (\sigma_J )
(adx_{j_1}\ldots dx_{j_{r}})(a'dx_{j_{r+1}}\ldots
dx_{j_n}).
\end{array}
\]
Hence, we have
\[
p_J(Q_{i_{n-r}}\ldots Q_{i_1} dx_1\ldots dx_n)
=
(-1)^{(n-r)r} \mathrm{sgn}(\sigma_J) Q_{i_{n-r}}\ldots Q_{i_1} 
dx_{j_{r+1}}\ldots dx_{j_n}
\]
as required.
\end{proof}

\bibliographystyle{gtart}
\bibliography{link}

\begin{thebibliography}{}
\providecommand\bibmarginpar{\leavevmode\marginpar}
\def\urlstyle#1{{\tt #1}}

\bibitem{andersen}
\textbf{K Andersen}, \textbf{J Grodal}, \textbf{J M{\o}ller}, \textbf{A
  Viruel}, \emph{The classification of $p$--compact groups for $p$ odd}
  \xox{arXiv}{math.AT/0302346}

\bibitem{benson}
\textbf{D\,J Benson}, \emph{Polynomial invariants of finite groups}, London
  Mathematical Society Lecture Note Series 190, Cambridge University Press,
  Cambridge (1993) \xox{MR}{1249931}

\bibitem{crabb}
\textbf{M\,C Crabb}, \href{http://dx.doi.org/10.1112/S0024609305004704}
  {\emph{Dickson--{M}\`ui invariants}}, Bull. London Math. Soc. 37 (2005)
  846--856 \xox{MR}{2186717}

\bibitem{dickson}
\textbf{L\,E Dickson},
  \href{http://links.jstor.org/sici?sici=0002-9947(191101)12:1%3C75:AFSOIO%3E2%
.0.CO%3B2-%23} {\emph{A fundamental system of invariants of the general modular
  linear group with a solution of the form problem}}, Trans. Amer. Math. Soc.
  12 (1911) 75--98 \xox{MR}{1500882} \xox{JFM}{42.0136.01}

\bibitem{greiss}
\textbf{R\,L Griess, Jr}, \href{http://dx.doi.org/10.1007/BF00150757}
  {\emph{Elementary abelian {$p$}--subgroups of algebraic groups}}, Geom.
  Dedicata 39 (1991) 253--305 \xox{MR}{1123145}

\bibitem{kameko}
\textbf{M Kameko}, \textbf{M Mimura},
  \href{http://dx.doi.org/10.2140/gtm.2006.10.209} {\emph{On the
  Rothenberg-Steenrod spectral sequence for the $\mathrm{mod}~3$ cohomology of
  the classifying space of the exceptional Lie group $E_8$}}, from:
  ``Proceedings of the Nishida Fest (Kinosaki 2003)'', Geom. Topol. Monogr. 10
  209--222

\bibitem{mui}
\textbf{H M\`ui}, \emph{Modular invariant theory and cohomology algebras of
  symmetric groups}, J. Fac. Sci. Univ. Tokyo Sect. IA Math. 22 (1975) 319--369
  \xox{MR}{0422451}

\bibitem{smith}
\textbf{L Smith}, \emph{Polynomial invariants of finite groups}, Research Notes
  in Mathematics 6, A K Peters Ltd., Wellesley, MA (1995) \xox{MR}{1328644}

\bibitem{wilkerson}
\textbf{C Wilkerson}, \emph{A primer on the {D}ickson invariants}, from:
  ``Proceedings of the Northwestern Homotopy Theory Conference (Evanston, Ill.,
  1982)'', Contemp. Math. 19, Amer. Math. Soc., Providence, RI (1983)  421--434
  \xox{MR}{711066}\ \ A corrected version is available at the Hopf Topology
  Archive.

\end{thebibliography}

\end{document}